\documentstyle[amstex,12pt,leqno]{article}

\title{Painlev\'e equations and deformations of rational surfaces
with rational double points\footnote{1991 {\it Mathematical Subject Classification.} 14D15, 34M55, 32G10, 14J26}
\footnote{{\it Key Words and phrases.} Painlev\'e 
equations, 
B\"ackglund transformations, Rational Surfaces, Rational double points, Deformations} 
}
\author{Masa-Hiko SAITO
\footnote{Partly supported by Grant-in Aid
for Scientific Research (B-09440015), (B-12440008) and (C-11874008), 
the Ministry of
Education, Science and Culture, Japan }
\\
Department of Mathematics, Faculty of Science, Kobe 
University \\
\\
Hiroshi UMEMURA\footnote{Partly supported by Grant-in Aid
for Scientific Research (B-11440006), 
the Ministry of
Education, Science and Culture, Japan } \\
Graduate School of Mathematics, Nagoya University} 
\date{}

\begin{document}

\newcommand{\A}{{\mathbb A}}
\newcommand{\x}{{\mathfrak X}}
\newcommand{\C}{{\mathbb C}}
\renewcommand{\P}{{\mathbb P}}
\newcommand{\D}{{\mathfrak D}}
\newcommand{\Z}{{\mathfrak Z}}
\newcommand{\m}{{\mathfrak m}}
\newcommand{\Y}{{\cal Y}}
\newcommand{\Pic}{{\rm Pic} \ }
\newcommand{\Proof}{{\bf Proof}.\ \ }
\newcommand{\Spec}{{\rm Spec} \  }
\renewcommand{\mod}{{\rm mod}\ }
\newcommand{\Ext}{{\rm Ext}}

\renewcommand{\theequation}{\thenewsection.\arabic{equation}}

\newtheorem{Theorem}{\ \ \ Theorem}
\newtheorem{Remark}[Theorem]{\ \ \ Remark}
\newtheorem{Proposition}[Theorem]{\ \ \ Proposition}
\newtheorem{Question}[Theorem]{\ \ \ Question}
\newtheorem{Corollary}[Theorem]{\ \ \ Corollary}
\newtheorem{Lemma}[Theorem]{\ \ \ Lemma}
\newtheorem{Sublemma}[Theorem]{\ \ \ Sublemma}
\newcounter{meq}
\newcounter{newsection}
\renewcommand{\theTheorem}{\addtocounter{equation}{1} (\thenewsection.\arabic{equation}).}
\renewcommand{\thesection}{\S \arabic{section}.}
\renewcommand
{\thenewsection}{\setcounter{newsection}{\value{section}}\arabic{newsection}}
\maketitle

In this paper we give an answer to the fundamental questions about
the Painlev\'e equations.
Where do the B\"acklund transformations come from? 
Our approach depends on the geometry of the projective surface constructed
by Okamoto and reviewed in \cite{Utr}.

The Painlev\'e equations were discovered around 1900 in the pursuit 
of special functions.
Painlev\'e and Gambier classified algebraic differential equations
$y''=R(t,y,y')$ without movable  singular points, where $R$ is a rational
function of $t,y,y'$ and $t$ is the independent variable.
After the refinement procedure of throwing away 
among 
the differential equations those that are integrable by the so far known functions,
they arrived at the list of the six Painlev\'e equations $P_J$ 
$(I \leq J \leq VI )$.
After the discovery of Painlev\'e equations,
an unexpected but important feature of the Painlev\'e equations was revealed.
The Painlev\'e equations have symmetries or the B\"acklund transformations.
In fact, if we recall the motivation of the discovery, it is surprising that
the Painlev\'e equations admit symmetries.
So a natural question arise.
Where do the B\"acklund transformations come from?
Our answer is that they arise from the rational double points.

The Painlev\'e equation $P_J$ is written in a Hamiltonian form
\begin{eqnarray*}
\left\{
\begin{array}{lcr}
\frac{\displaystyle dq}{\displaystyle dt} &=& 
\frac{\displaystyle \partial H_J}{\displaystyle \partial p},\\
\smallskip \\
\frac{\displaystyle dp}{\displaystyle dt} &=& 
-\frac{\displaystyle \partial H_J}{\displaystyle \partial q}
\end{array}
\right.
\end{eqnarray*}
for an appropriate polynomial $H_J$ of $q$, $p$ and $t$ $(I \leq J \leq VI)$.
We also clarify where the polynomial $H_J$ comes from and 
why it is  uniquely determined.
This point is a geometric interpretation of \cite{ST} and \cite{MMT}.  
It is interesting to notice that our theory provides us a
lot of pairs $(U, \, V)$ of algebraic varieties defined over $\C$ such that
$U$ and $V$ are not isomorphic as algebraic varieties but they are isomorphic 
 as complex manifolds (cf \S 11).  
 (The idea in \S 11 
is developed further in the direction of deformation theory of Okamoto--Painlev\'e pairs and its relation to Painlev\'e equations in \cite{STT}.)

\par
We work with the second Painlev\'e equation $P_{II}$ to 
illustrate
the general case.
Our argument is given in a form easily applicable to the other Painlev\'e equations (cf. \cite{STT}).  
Through out the paper, the ground field in $\C$.
The most natural setting seems to be, however, 
over the ring ${\mathbb Z}[\frac{1}{2}]$. 

We should mention that a recent work of Sakai 
\cite{Sakai} is also  working  on  the discrete and differential Painlev\'e equations related to the geometry of rational surfaces and symmetries of affine Weyl groups.  

The authors express gratitude to Professor Kei-ichi Watanabe. The discussions with him was indispensable to write \S 9.

\section{Construction of a family of rational surfaces.}

\par
The second Painlev\'e equation is an ordinary differential equation
\begin{eqnarray*}
P_{II}(\alpha)  \hspace{1cm}  y''=2y^3+ty+\alpha
\end{eqnarray*}
of the second order,
where $t$ is the independent variable, $y''=d^2y/dt^2$ and
$\alpha \in \C$ is a parameter.
The extended affine Weyl group $\widehat{G}$ of type
$A^{(1)}_1$ appears as the symmetry of $P_{II}(\alpha)$.
Namely if y is a solution of $P_{II}(\alpha)$, then 
\begin{eqnarray*}
T_+(y)=-y-\frac{\alpha+\frac{1}{2}}{y'+y^2+\frac{t}{2}}
\end{eqnarray*}
is a solution of $P_{II}(\alpha+1)$,
\begin{eqnarray*}
T_-(y)=-y+\frac{\alpha-\frac{1}{2}}{y'-y^2-\frac{t}{2}}
\end{eqnarray*}
is a solution of $P_{II}(\alpha-1)$ and 
\begin{eqnarray*}
I(y)=-y
\end{eqnarray*}
is a solution of $P_{II}(-\alpha)$.
The automorphisms $i$, $t_+$, $t_-$ of the affine line $\A^1$ with
coordinate system $\alpha$, i.e. $\A^1=\Spec \C[\alpha]$,
given by
\begin{eqnarray*}
t_+(\alpha)=\alpha+1, \hspace{0.5cm} t_-(\alpha)=\alpha-1, \hspace{0.5cm}
i(\alpha)=-\alpha
\end{eqnarray*} 
generate the (extended) affine Weyl group $\widehat{G}$ of type
$A^{(1)}_1$.
So the group $\widehat{G}$ operates on the affine line $\A^1$.
We extend the operation of $\widehat{G}$ on the affine plane $\A^2$
with coordinate system $(t, \alpha)$ so that $\widehat{G}$ leaves $t$ 
invariant.
We consider the affine space $\A^4$ with coordinate system 
$(t, \alpha, y ,y')$ and the projection 
$p_{12}:\A^4 \rightarrow \A^2, 
\ (t, \alpha,y,y') \mapsto (t,\alpha)$. 
Through the transformations $T_+$, $T_-$, $I$ regarded as birational 
automorphism of the affine space $\A^4$, the extended affine Weyl group
$\widehat{G}$ operates birationally on $\A^4$
such that the projection $p_{12}:\A^4 \rightarrow \A^2$
is $\widehat{G}$-equivariant.
For example the birational automorphism of $\A^4$
induced by $T_+$ is given by
\begin{eqnarray*}
\A^4 &\rightarrow&\A^4,\\
\smallskip \nonumber \\
(t,\alpha, y, y') &\mapsto& (t, \alpha + 1, T_+(y), T_+(y)' )
\end{eqnarray*}
where
\begin{eqnarray*}
T_+(y)=-y-\frac{\alpha+\frac{1}{2}}{y'+y^2+\frac{t}{2}} 
\end{eqnarray*}
as is given above and
\begin{eqnarray*}
T_+(y)'&=&\left(-y-\frac{\alpha+\frac{1}{2}}{y'+y^2+\frac{t}{2}}\right)' 
\\
\smallskip \nonumber\\
&=&-y'+\frac{(\alpha+\frac{1}{2})(y''+2yy'+\frac{1}{2})}
{(y'+y^2+\frac{t}{2})^2}\\
\smallskip \nonumber \\ 
&=&-y'+\frac{(\alpha+\frac{1}{2})(2y^3+ty+\alpha+2yy'+\frac{1}{2})}
{(y'+y^2+\frac{t}{2})^2}.
\end{eqnarray*}
Namely if we use a more rigorous notation, 
the birational automorphism of $\A^4$ is defined by 
$$
\A^4 \rightarrow \A^4,\qquad 
(t,\alpha + 1, u, v) \mapsto (t, -\alpha, U(t,\alpha, u,v), V(t,\alpha, u, v))
$$
where
\begin{eqnarray*}
U(t,\alpha, u,v)&=&-u-\frac{\alpha+\frac{1}{2}}{v+u^2+\frac{t}{2}},\\
\smallskip \nonumber\\
V(t, \alpha, u, v) 
&=&-v+\frac{(\alpha+\frac{1}{2})(2u^3+tu+\alpha+2uv+\frac{1}{2})}
{(v+u^2+\frac{t}{2})^2}.
\end{eqnarray*}
\begin{Remark}
For the root system of type {\rm $A_1$}, 
the affine Weyl group is isomorphic to the extended affine 
Weyl group. 
For a systematic understanding of the Painlev\'e equations,
we should regard {\rm $\widehat{G}$} as the extended affine Weyl group. 
\end{Remark}

On the affine space $\A^4$, we have a vector field
\begin{eqnarray*}
\delta(\alpha)=\frac{\partial }{\partial t}
+y'\frac{\partial}{\partial y}
+(2y^3+ty+\alpha)\frac{\partial}{\partial y'}
\end{eqnarray*}
so that the birational maps $T_+$, $T_-$, $I$ are compatible
with the vector field $\delta(\alpha)$.
We projectified the fibration $p=p_{12}:\A^4 \rightarrow \A^2$
in the following manner (\cite{Utr}).
To explain the projectification of $p$, we had better take other
coordinate systems on $\A^4$ and $\A^2$.
First, we introduce a new coordinate system $(t,c)$ such that $c=\alpha-\frac{1}{2}$. 
Second, we define a new coordinate system $(t,c,p,q)$ on $\A^4$
such that 
\begin{eqnarray}
\left\{
\begin{array}{lcl}
q &=& y,\\
\smallskip \\
p &=& y'-q^2-\frac{t}{2}.
\end{array}
\right.
\end{eqnarray} 
In terms of this coordinate system, the second Painlev\'e equation
is written as
\begin{eqnarray*}
S_2(c)
\left\{
\begin{array}{lcl}
\frac{
\displaystyle 
dq}{
\displaystyle
dt} &=& 
q^2+p+\frac{ 
\displaystyle 
t}{
\displaystyle 
2},\\
\smallskip \\
\frac{
\displaystyle 
dp}{
\displaystyle 
dt} &=& -2qp+c
\end{array}
\right.
\end{eqnarray*}
and the vector field $\delta(\alpha)$ is given by
\begin{eqnarray*}
D(c)=\frac{\partial}{\partial t}
+(q^2+p+\frac{t}{2})\frac{\partial}{\partial q}
+(-2qp+c)\frac{\partial}{\partial p}.
\end{eqnarray*}
In other words, if we consider an isomorphism
\begin{eqnarray*}
\varphi:\A^4 &\rightarrow& \A^4\\
\smallskip \nonumber\\
(t,\alpha, y, y') &\mapsto& \left(t,c,y,y'-y^2-\frac{t}{2}\right),
\end{eqnarray*}
then $\varphi$ transforms $\delta(\alpha)$ to $D(c)$.
The system $S_2(c)$ is a Hamiltonian form. In fact if 
we take $H(t,c,q,p)=q^2p+\frac{1}{2}p^2+\frac{t}{2}p-cq$,
then we have 
\begin{eqnarray*}
S_2(c)
\left\{
\begin{array}{lcl}
\frac{\displaystyle dq}{\displaystyle dt} 
&=&
\frac{\displaystyle \partial H}{\displaystyle \partial p}\\
\smallskip \\
\frac{\displaystyle dp}{\displaystyle dt} 
&=&
-\frac{\displaystyle \partial H}{\displaystyle \partial q}.
\end{array}
\right.
\end{eqnarray*}
We explain now how to projectify the fibration 
$p_{12}:\A^4 \rightarrow \A^2$.
We take a point $(t_0,c_0) \in \A^2$ fixed once for all
and show how to projectify the fiber 
$X[t_0, c_0]:=p_{12}^{-1}(t_0,c_0)$ of the morphism
$p_{12}:=\A^4 \rightarrow \A^2$, $(t,c,q,p) \mapsto (t,c)$.
The fiber $X[t_0, c_0]$ is isomorphic to the affine plane $\A^2$.
We need four copies $W_i$ $(1 \leq i \leq 4)$ of $\A^2$ 
with coordinate system $(y_i,z_i)$.
We glue together the $W_i$' s by the following rule to get a rational ruled
surface $Z[t_0, c_0]$.
\\{\bf (i)} 
A point $(y_1,z_1) \in W_1$ and a point $(y_2,z_2) \in W_2$ are
identified if 
\begin{eqnarray*}
y_1=y_2 \hspace{0.5cm} \mbox{and} \hspace{0.5cm} z_1z_2=1.
\end{eqnarray*}
\\{\bf(ii)} 
A point $(y_1,z_1) \in W_1$ and a point $(y_3,z_3) \in W_3$ are
identified if 
\begin{eqnarray*}
y_1 y_3 =1 \hspace{0.5cm} \mbox{and} \hspace{0.5cm} z_1 = c_0 y_3-y_3^2z_3.
\end{eqnarray*}
We notice that the latter condition is equivalent to
\begin{eqnarray*}
z_3=c_0y_1-y_1^2 z_1.
\end{eqnarray*}
\\{\bf(iii)} 
A point $(y_3,z_3) \in W_3$ and a point $(y_4,z_4) \in W_4$ are
identified if 
\begin{eqnarray*}
y_3=y_4 \hspace{0.5cm} \mbox{and} \hspace{0.5cm} z_3z_4=1.
\end{eqnarray*}

The projections 
\begin{eqnarray*}
W_i \rightarrow \A^1, \hspace{0.5cm} (y_i, z_i) \mapsto y_i
\hspace{0.5cm} \mbox{for $1 \leq i \leq 4$}
\end{eqnarray*}
glue together to give a fibration
\begin{eqnarray*}
p_{Z[t_0,c_0]}:Z[t_0,c_0] \rightarrow \P^1.
\end{eqnarray*}
So $Z[t_0,c_0]$ is a $\P^1$-bundle over $\P^1$ 
on 
$Z[t_0,c_0]$ is a rational ruled surface.
We have, as is easily seen 
\begin{eqnarray*}
Z[t_0,c_0]\cong \left\{ 
\begin{array}{ll}
F_2 \hspace{0.5cm}&\mbox{if $c_0=0$},\\
\smallskip \\
\P^1 \times \P^1 &\mbox{otherwise}.
\end{array}
\right.
\end{eqnarray*}

It is apparent that the construction of $Z[t_0,c_0]$ depends only on $c_0$
and not on $t_0$.
We identify the fiber $X[t_0,c_0]$, which is the affine plane
with coordinate system $(q,p)$, with $W_1$ that is the affine space 
with coordinate system $(y_1, z_1)$ by sending a point 
$(q,p)$ of $X[t_0, c_0]$ to the point $(y_1,z_1)=(q,p)$ of $W_1$.
We thus constructed a projectification $Z[t_0, c_0]$ of the fiber $X[t_0,c_0]$ 
but this is not the desired projectification.
To get the projectification $X[t_0,c_0]$, we have to blow-up 
the rational ruled surface $Z[t_0,c_0]$ eight times, the centers being carefully
chosen infinitely near points of $(y_4, z_4)=(0,0) \in W_4$.
The centers depend not only on $c_0$ but also $t_0$ so that the 
projectification $\x[t_0, c_0]$ depends on $c_0$ and $t_0$.
The center $a_1$ of the first blowing-up is the point $(y_4, z_4)=(0,0)$
on $W_4$.
For $1 \leq i \leq 8$, we denote by 
$\pi_i:Z_i[t_0,c_0] \rightarrow Z_{i-1}[t_0,c_0]$ the $i$-the
blowing-up, by $E_i \subset Z_i[t_0,c_0]$ the exceptional divisor of $\pi_i$,
and by $a_i \in Z_{i-1}[t_0,c_0]$ the center of the $i$-th blown-up $\pi_i$,
so that $Z_0[t_0,c_0]=Z[t_0,c_0]$, $Z_8[t_0,c_0]=\x[t_0,c_0]$ 
and $a_1=(0,0) \in W_4$.  
The curve $S=\{(y_2,z_2) \in W_3 \, | \,   z_3=0\} \cup \{(y_4,z_4)\, | \,  z_4=0\}$
is a section of the ruled surface $\pi:Z[t_0,z_0] \rightarrow \P^1$
and the self-intersection number $S^2=2$.
For $2 \leq i \leq 8 $ the center $a_i \in Z_{i-1}[t_0,c_0]$ is always 
on the exceptional divisor $E_{i-1}$ of the previous blowing-up.
For $i=2,3,4$, the center $a_i\in Z_{i-1}[t_0,c_0]$ is the intersection
point of $E_{i-1}$ and the proper transform of $S$ by the morphism 
$Z_{i-1}[t_0,c_0] \rightarrow Z_0[t_0,c_0]=Z[t_0,c_0]$.
Let $D_i \subset \x[t_0,c_0]=Z_8[t_0,c_0]$ be proper transform of $E_i$
by the morphism $\x[t_0, c_0]=Z_8[t_0,c_0] \rightarrow Z_i[t_0,c_0]$
for $1 \leq i \leq 7$.
The following result is proved in \cite{Utr}.
\vskip 5mm
{\bf Sublemma }(3.6) {\bf of} \cite{Utr}.
{\it Locally on $W_4$, the construction 
$\x[t_0,c_0]=Z_8[t_0,c_0] \rightarrow Z[t_0,c_0]$ or the blowing-up
of $Z[t_0,c_0]$ is equivalent to the minimal resolution of the rational map
\begin{eqnarray*}
F:W_4 \cdots \rightarrow \P^1, \hspace{0.5cm}
(y_4,z_4) \mapsto (y_4^4z_4, 2z_4-y_4^4+ty_4^2z_4+(2c_0+1)y_4^3z_4).
\end{eqnarray*}
}

{\bf Remark.}
{\it $F$ is not regular on $\x[t_0,c_0]$.
For, it has a base point $(y_2,z_2)=0$ on $W_2$.
This point outside $W_4$ is left untouched in the construction of 
$\x[t_0, c_0]$.
}

\vskip 5mm
We denote by $D_0$ the proper transform of the curve $S \subset Z[t_0, c_0]$
by the morphism $\x[t_0,c_0] \rightarrow Z[t_0,c_0]$.
Then we have $D_1 \simeq \P^1$ and $D_i^2=-2$ for $ 1\leq i\leq 7 $.
The configuration of the $D_i$'s is shown in Fig.(1.3).

\unitlength 1.5mm
\begin{figure}[h]
\begin{center}
\begin{picture}(65,33)

\put(6,2){$D_1$}
\put(0,15){$D_2$}
\put(16,10){$D_3$}
\put(12,25){$D_4$}
\put(30,11){$D_0$}
\put(44,11){$D_5$}
\put(62,15){$D_6$}
\put(52,2){$D_7$}

\put(8,5){\line(0,1){14}}
\put(6,16){\line(1,0){14}}
\put(18,14){\line(0,1){14}}
\put(16,26){\line(1,0){32}}
\put(32,14){\line(0,1){14}}
\put(46,14){\line(0,1){14}}
\put(44,16){\line(1,0){14}}
\put(56,5){\line(0,1){14}}
\end{picture}\\
{\rm Fig.(1.3)}
\end{center} 
\end{figure}

The dual graph of Fig.(1.3) is Fig.(1.4) that is the Dynkin diagram of type 
$E_7^{(1)}$.
\begin{figure}[h]
\unitlength 1.5mm
\begin{center}
\begin{picture}(70,18)
\multiput(5,15)(10,0){7}{\circle{1}}
\multiput(5.5,15)(10,0){6}{\line(1,0){9}}
\put(35,14.5){\line(0,-1){9}}
\put(35,5){\circle{1}}
\put(3,17){$D_1$}
\put(13,17){$D_2$}
\put(23,17){$D_3$}
\put(33,17){$D_4$}
\put(43,17){$D_5$}
\put(53,17){$D_6$}
\put(63,17){$D_7$}
\put(36,5){$D_0$}
\end{picture}\\
{\rm Fig.(1.4)}
\end{center}
\end{figure}
\addtocounter{equation}{2}

Namely the vertices in Fig.(1.4) correspond to the curves in Fig.(1.3).
Two different vertices are jointed by an edge if corresponding two curves have a point in common.
We set $\D[t_0, c_0]:=\cup_{i=0}^7D_i \subset \x[t_0,c_0]$, which is a divisor
on $\x[t_0,c_0]$.
The divisors $D_i$ also depend on $t_0$ and $c_0$.
So we should denoted $D_i$ by $\D_i[t_0,c_0]$.
To simplify the notation, we use $\D_i[t_0,c_0]$, only when we emphasize
the dependence on $t_0$ and $c_0$.
So far we fixed $(t,c)=(t_0,c_0)$.
The above construction works globally on the fibration 
$\pi:\A^4 \rightarrow \A^2$
and gives us a fiber space
\begin{eqnarray*}
\varphi:\x \rightarrow \A^2
\end{eqnarray*}
and a divisor $\D$ on $\x$ such that for a point 
$(t_0,c_0) \in \A^2$, the fiber $\varphi^{-1}(t_0,c_0)$ is isomorphic to
$\x[t_0,c_0]$ and such that $\D \cap \varphi^{-1}(t_0,c_0)$
yields the divisor $\D[t_0,c_0]$ on $\x[t_0,c_0]$.
Namely, let $\psi$ be the composite $p_2 \circ \varphi$
\begin{eqnarray*}
\x \stackrel{\varphi}{\rightarrow} \A^2 \stackrel{p_2}{\rightarrow} \A,
\end{eqnarray*}
where $p_2:\A^2 \rightarrow \A^1$ is the projection onto the second 
factor so that $p_2(t,c)=c$ for a point $(t,c) \in \A^2$.
We denote the fiber $\psi^{-1}(c_0)$ by $\x[c_0]$ 
for a point $c_0 \in \A^1$.
Then $\varphi$ gives a morphism
\begin{eqnarray*}
\varphi_{c_0}:\x[c_0] \rightarrow \A^1\times c_0 \subset  \A^2.
\end{eqnarray*}
Since $\A^1 \times c_0 \simeq \Spec \A[t]$
is the affine line with coordinate system $t$,
we have a morphism
\begin{eqnarray*}
\varphi[c_0]:\x[c_0] \rightarrow \Spec \C[t],
\end{eqnarray*}
which we denote also by $\varphi_{c_0}$.
On the threefold $\x[c_0]$, we have a divisor 
$\D[c_0]:=\D \cap \x[c_0]$ 
and a rational vector field
\begin{eqnarray*}
\delta[c_0]=\frac{\partial}{\partial t}
+(y_1^2+z_1+\frac{t}{2})\frac{\partial}{\partial y_1}
+(-2 y_1 z_1+c_0)\frac{\partial}{\partial z_1}
\end{eqnarray*}
generates a foliation on $\x[c_0]$.
The divisor $\D[c_0]$ is the set of singular points of the foliation so that
every leaf is transversal to the fibers of $\x \rightarrow \Spec \C[t]$
in the open set $\x[c_0] \backslash \D[c_0]$.
Since the Painlev\'e equations have no movable singular points
(=movable branch points and movable essential singular points),
the set of solutions on the set of leaves of the foliation in 
$\x[c_0] \backslash \D[c_0] $ sweeps
the whole space $\x[c_0] \backslash \D[c_0]$ (See Fig.(1.5)).

\begin{figure}[h]

\begin{center}
\unitlength 0.1in
\begin{picture}(30.00,30.50)(6.20,-30.60)
%
\special{pn 8}%
\special{ar 1140 1270 520 1090  1.5563465 6.2831853}%
\special{ar 1140 1270 520 1090  0.0000000 1.5058472}%
%
\special{pn 8}%
\special{ar 3400 1200 520 1090  4.5975124 6.2831853}%
\special{ar 3400 1200 520 1090  0.0000000 1.5058472}%
%
\special{pn 8}%
\special{pa 1080 350}%
\special{pa 1090 710}%
\special{fp}%
\special{pa 980 530}%
\special{pa 1240 510}%
\special{fp}%
\special{pa 1170 420}%
\special{pa 1180 720}%
\special{fp}%
\special{pa 970 440}%
\special{pa 960 730}%
\special{fp}%
\special{pa 920 520}%
\special{pa 1050 520}%
\special{fp}%
\special{pa 860 640}%
\special{pa 1000 660}%
\special{fp}%
\special{pa 1140 670}%
\special{pa 1380 620}%
\special{fp}%
\special{pa 1310 560}%
\special{pa 1320 790}%
\special{fp}%
\special{pa 870 560}%
\special{pa 870 850}%
\special{fp}%
%
\special{pn 8}%
\special{pa 3360 290}%
\special{pa 3370 650}%
\special{fp}%
\special{pa 3260 470}%
\special{pa 3520 450}%
\special{fp}%
\special{pa 3450 360}%
\special{pa 3460 660}%
\special{fp}%
\special{pa 3250 380}%
\special{pa 3240 670}%
\special{fp}%
\special{pa 3200 460}%
\special{pa 3330 460}%
\special{fp}%
\special{pa 3140 580}%
\special{pa 3280 600}%
\special{fp}%
\special{pa 3420 610}%
\special{pa 3660 560}%
\special{fp}%
\special{pa 3590 500}%
\special{pa 3600 730}%
\special{fp}%
\special{pa 3150 500}%
\special{pa 3150 790}%
\special{fp}%
%
\special{pn 8}%
\special{ar 1870 1400 500 420  3.3665649 5.3381592}%
%
\special{pn 8}%
\special{ar 3070 1000 960 180  1.0776007 2.6831437}%
%
\special{pn 8}%
\special{pa 3420 1170}%
\special{pa 3630 1140}%
\special{fp}%
\special{sh 1}%
\special{pa 3630 1140}%
\special{pa 3561 1130}%
\special{pa 3577 1148}%
\special{pa 3567 1169}%
\special{pa 3630 1140}%
\special{fp}%
%
\special{pn 8}%
\special{ar 1787 1780 500 420  3.3665649 5.3381592}%
%
\special{pn 8}%
\special{ar 2987 1380 960 180  1.0776007 2.6831437}%
%
\special{pn 8}%
\special{pa 1110 180}%
\special{pa 3400 110}%
\special{fp}%
\special{pa 1150 2340}%
\special{pa 3440 2290}%
\special{fp}%
%
\special{pn 8}%
\special{pa 3310 1560}%
\special{pa 3580 1510}%
\special{fp}%
\special{sh 1}%
\special{pa 3580 1510}%
\special{pa 3511 1502}%
\special{pa 3528 1520}%
\special{pa 3518 1542}%
\special{pa 3580 1510}%
\special{fp}%
%
\special{pn 8}%
\special{pa 3460 2290}%
\special{pa 3428 2295}%
\special{pa 3396 2297}%
\special{pa 3364 2296}%
\special{pa 3333 2289}%
\special{pa 3303 2279}%
\special{pa 3274 2266}%
\special{pa 3246 2251}%
\special{pa 3219 2232}%
\special{pa 3195 2211}%
\special{pa 3172 2189}%
\special{pa 3149 2166}%
\special{pa 3128 2142}%
\special{pa 3109 2117}%
\special{pa 3091 2090}%
\special{pa 3075 2063}%
\special{pa 3058 2035}%
\special{pa 3042 2007}%
\special{pa 3028 1979}%
\special{pa 3015 1950}%
\special{pa 3002 1920}%
\special{pa 2990 1891}%
\special{pa 2979 1861}%
\special{pa 2968 1831}%
\special{pa 2958 1800}%
\special{pa 2949 1770}%
\special{pa 2940 1739}%
\special{pa 2932 1708}%
\special{pa 2924 1677}%
\special{pa 2917 1646}%
\special{pa 2910 1614}%
\special{pa 2904 1583}%
\special{pa 2899 1551}%
\special{pa 2894 1520}%
\special{pa 2890 1488}%
\special{pa 2886 1456}%
\special{pa 2882 1424}%
\special{pa 2879 1393}%
\special{pa 2876 1361}%
\special{pa 2874 1329}%
\special{pa 2872 1297}%
\special{pa 2871 1265}%
\special{pa 2870 1233}%
\special{pa 2870 1201}%
\special{pa 2870 1169}%
\special{pa 2871 1137}%
\special{pa 2872 1105}%
\special{pa 2874 1073}%
\special{pa 2876 1041}%
\special{pa 2878 1009}%
\special{pa 2881 977}%
\special{pa 2884 945}%
\special{pa 2887 914}%
\special{pa 2892 882}%
\special{pa 2897 850}%
\special{pa 2902 819}%
\special{pa 2908 787}%
\special{pa 2914 756}%
\special{pa 2921 725}%
\special{pa 2928 693}%
\special{pa 2936 662}%
\special{pa 2944 631}%
\special{pa 2954 601}%
\special{pa 2964 570}%
\special{pa 2974 540}%
\special{pa 2986 510}%
\special{pa 2998 481}%
\special{pa 3011 451}%
\special{pa 3024 423}%
\special{pa 3039 394}%
\special{pa 3054 366}%
\special{pa 3071 338}%
\special{pa 3088 311}%
\special{pa 3107 285}%
\special{pa 3126 260}%
\special{pa 3148 236}%
\special{pa 3170 214}%
\special{pa 3194 192}%
\special{pa 3219 173}%
\special{pa 3247 156}%
\special{pa 3275 141}%
\special{pa 3305 129}%
\special{pa 3336 122}%
\special{pa 3347 120}%
\special{sp -0.045}%
%
\special{pn 8}%
\special{pa 2160 1050}%
\special{pa 2200 1070}%
\special{dt 0.045}%
\special{pa 2200 1070}%
\special{pa 2199 1070}%
\special{dt 0.045}%
\special{pa 2170 1060}%
\special{pa 2300 1110}%
\special{dt 0.045}%
\special{pa 2300 1110}%
\special{pa 2299 1110}%
\special{dt 0.045}%
\special{pa 2040 1420}%
\special{pa 2160 1480}%
\special{dt 0.045}%
\special{pa 2160 1480}%
\special{pa 2159 1480}%
\special{dt 0.045}%
\put(0,-28){\line(1,0){40}}
\put(41,-29){$\Spec \C[t]$}
\put(35,-24){\vector(0,-1){2}}
\put(36,-27){$\varphi_{c_0}$}
\put(41,-5){$\x[c_0]$}
\end{picture}\\
Fig.(1.5)
\end{center}
\end{figure}
\addtocounter{equation}{1}
$\x[t_0,c_0] \backslash \D[t_0,c_0]$ is the space of initial conditions
of the second Painlev\'e equation introduced by Okamoto \cite{Ofol}
(cf. \cite{Utr}).
In terms of coordinate system $(t,c,q,p)$ of the affine space $\A^4$,
the affine extended Weyl group $G$ is generated by rational maps
\begin{eqnarray*}
J:\A^4 \cdots \rightarrow \A^4, \ \ (t,c,q,p) \mapsto (t,-1-c,-q,-2q^2-p-t),
\end{eqnarray*}
and 
\begin{eqnarray*}
I:\A^4 \cdots \rightarrow \A^4, \ \  
(t,c,q,p) \mapsto (t,\,-c,\, q-\frac{c}{p},\, p).
\end{eqnarray*}
The corresponding birational map induced respectively by $J$
(resp. $I$) on a  variety $\Z$ over $\A^2$, which
is $\A^2$-birational to ${p_{12}}:\A^4 \rightarrow \A^2$,
will be denoted by $J_{\Z}(c, -1-c)$ (resp. $I_{\Z}(c, -c)$).
In this context, when we specialize the parameter $c$ to $c_0 \in \C$,
we denote the thus obtained birational map by $J_{\Z}(c_0, -1-c_0)$,
$I_{\Z}(c_0,-c_0)$.
We proved in \cite{Utr} that when the parameter $c$ takes a fixed value 
$c_0 \in \C$,
the B\"acklund transformations give the following 
isomorphisms that commute with the foliations.
\begin{eqnarray}
J_{\x}(c_0,-1-c_0):\x[c_0] \rightarrow \x[-1-c_0].\\
\smallskip \nonumber\\
I_{\x}(c_0,-c_0):\x[c_0] \rightarrow \x[-c_0].
\end{eqnarray}
Moreover we have 
\begin{eqnarray}
J_{\x}(c_0,-1-c_0)[\D[c_0]]=\D[-1-c_0]
\end{eqnarray}
and 
\begin{eqnarray}
I_{\x}(c_0,-c_0)[\D[c_0]]=\D[-c_0].
\end{eqnarray}
For, for every $c_0 \in \C$, $\D[c_0]$ 
is the set of singular points of the foliation on $\x[c_0]$
and the isomorphisms $J_{\x}(c_0,-1-c_0)$ and $I_{\x}(c_0,-c_0)$
commute with the foliation.
\begin{Proposition}
We have 
\begin{eqnarray*}
J_{\x}(c_0,-1-c_0)(\D_i[c_0])&=&\D_{8-i}[-1-c_0] 
\hspace{0.5cm} \mbox{for $1 \leq i \leq 7$},\\
\smallskip \nonumber\\
J_{\x}(c_0,-1-c_0)(\D_0[c_0])&=&\D_{0}[-1-c_0],\\
\smallskip \nonumber\\
I_{\x}(c_0,-c_0)(\D_i[c_0])&=&\D_i[-c_0]
\hspace{0.5cm} \mbox{for $0 \leq i \leq 7$}.
\end{eqnarray*}
\end{Proposition}

\Proof 
Once we have (1.8) and (1.9), we can check the identities easily.
\vspace{1pc}
\par
When we specialize not only $c$ to $c_0$ but also $t$ to $t_0$, 
the corresponding birational maps are denoted respectively by
$J_{\Z}(t_0;c_0,-1-c_0)$ and $I_{\Z}(t_0;c_0,-c_0)$.
Using this notation, we proved in \cite{Utr} that we have isomorphisms
\begin{eqnarray*}
J_{\x}(t_0;c_0,-1-c_0):\x[t_0,c_0] \rightarrow \x[t_0,-1-c_0]
\end{eqnarray*}
and 
\begin{eqnarray*}
I_{\x}(t_0;c_0,-c_0):\x[t_0,c_0] \rightarrow \x[t_0,-c_0].
\end{eqnarray*}
The following results follows from Proposition (1.10).

\begin{Proposition}
When $t$ and $c$ take fixed values $t_0$ and $c_0$, 
we have the following identities.
\begin{eqnarray*}
&&\hspace{-6cm}
{\bf (i)} \ \ J_{\x}(t_0;c_0,-1-c_0)(\D_i[t_0,c_0])=\D_{8-i}[t_0,-1-c_0]\\
\smallskip \\
&&\hspace{-6cm}\mbox{for $1 \leq i \leq 7$ and }\\
\smallskip \\
&&\hspace{-6cm}{\bf (ii)} \ \ J_{\x}(t_0;c_0,-1-c_0)(\D_0[t_0,c_0])=\D_0[t_0,-1-c_0].\\
\smallskip \\
&&\hspace{-6cm}{\bf(iii)} \ \ I_{\x}(t_0;c_0,-c_0)(\D_i)[t_0,-c_0].
\end{eqnarray*}
\end{Proposition}
\begin{Remark}
$J_{\x}; \x \rightarrow \x$ is a regular automorphism 
and the morphism $\varphi:\x \rightarrow \A^2$ is
compatible with automorphism $j:\A^2 \rightarrow \A^2,\ (t,c) \mapsto (t,-c)$. 
On the other hand, $I_{\x}:\x \rightarrow \x$ is not biregular.
\end{Remark}

In fact there exists a curve $C'_2[t_0,0]$ on $\x[t_0,0] \backslash \D$ 
isomorphic to $\P^1$ giving the Riccati solutions of 
$P_{II}\left(\frac{1}{2}\right)$.
The curve $C'_2[t_0,0]$ is the base locus of the rational map 
$I_{\x}:\x \rightarrow \x$. We proved in \cite{Utr}
that if we denote by $\widetilde{\x}$ the blowing-up of $\x$ along the surface 
$\bigcup_{t_0 \in \C} C'_2[t_0,0]$,
then $I_{\widetilde{\x}}:\widetilde{\x} \rightarrow \widetilde{\x}$ 
is a biregular automorphism.
We explain \S3 the reason why we have to blow up $\x$ to make $I_{\x}$
biregular.
The birational maps 
\begin{eqnarray*}
I_{\x}(c_0,-c_0):\x[c_0] \rightarrow \x[-c_0]
\end{eqnarray*}
and 
\begin{eqnarray*}
I_{\x}(t;c_0,-c_0):\x[t_0,c_0] \rightarrow \x[t_0,c_0]
\end{eqnarray*}
induced by $I_{\x}$ is the identity, when $c_0=0$.
So they are biregular.
Hence strictly speaking in the proofs of Propositions (1.10), (1.11)
we must treat the case $c_0=0$ separately.
For every integer $n$, there is curve $C[t_0,n]$
isomorphic to $\P^1$ on $\x[t_0,n]\backslash \D$. We blow up
$\x$ along the surfaces $\bigcup_{t_0\in \C}C[t_0,n] $ 
for every $n \in {\bf Z}$.
We get a manifold $\Y$, a projective limit of 
a scheme over $\A^2$ such that the extended Weyl group
$\widehat{G}=<I_{\Y},J_{\Y}>$ operates regularly on $\Y$.
In fact, we have seen in \cite{Utr},
$I_{\Y}$ and $J_{\Y}$ are automorphisms of $\Y$.
The extended affine Weyl group $\widehat{G}$ operates on $\Y$
and $\A^2$ in such a way that $\varphi_{\Y}:\Y \rightarrow \A^2$
is $\widehat{G}$-equivariant.
Here $\varphi_{\Y}:\Y \rightarrow \A^2$ is the composite of the blowing-up
morphism $\Y \rightarrow \x$ 
and the morphism $\varphi: \x \rightarrow \A^2$
arising from the projection $p_{12}: \A^4 \rightarrow \A^2$.
It is natural to ask a 
\begin{Question}
The quotient spaces $\Y/\widehat{G}$, $\A^2/\widehat{G}$
and the quotient morphism
$\Y/\widehat{G} \rightarrow \A^2/\widehat{G} $
are algebraizable?
\end{Question}

\section{Affine root systems on $\x[t_0,c_0]$.}
\setcounter{equation}{0}

\par
The additive group 
 $\Pic \x[t_0,c_0]$
 of linear equivalence classes
  of divisors on the surface 
$\x[t_0,c_0]$ is a lattice, of rank 10, the bilinear form on $\Pic \x[t_0,c_0]$
being defined by the intersection pairing.
\begin{Lemma}
We can blow down the surface $\x[t_0,c_0]$ to $\P^2$.
\end{Lemma}

\Proof 
If $c_0 \neq 0$, then the ruled surface $Z[t_0,c_0]$ is $\P^1 \times \P^1$.
Let $\widetilde{X}$ be the blown-up of $\P^1 \times \P^1$ 
at a point $P=(x_1,x_2)$ 
and let $C_1$, $C_2$ be respectively the proper transform of $x_1 \times \P^1$,
$\P^1 \times x_2$ under the blowing down morphism 
$\varphi:\widetilde{X} \rightarrow \P^1 \times \P^1$.
Then $C_1$ and $C_2$ are disjoint $-1$-curves so that we can collapse
$C_1$ and $C_2$ on $\widetilde{X}$ to get $\P^2$.
If the parameter $c_0 =0$, then the ruled surface $Z[t_0,c_0]$ is isomorphic
to $F_2$.
We blow-up $F_2$ at a point $P$ on the section $S$ of 
$\pi:F_2 \rightarrow \P^1$ with $S^2=2$ to get a surface $\widetilde{X}$.
So we have the blowing-down morphism $\varphi:\widetilde{X} \rightarrow F$.
Since $F_2 \backslash S=Z[t_0,c_0] \backslash S$ is a line bundle over $\P^1$,
we can find  a section $S_0$ of $\pi: F_2 \rightarrow \P^1$ such that $S_0^2=-2$and $(S.S_0)=0$.
Let $C_1$ be the proper transform of the fiber of the ruled surface 
$Z[t_0,c_0]$ passing through the point $P$ and let $C_2$ be
the proper transform of $S_0$ under $\varphi:\widetilde{X} \rightarrow F_2$.
Then since $C_1^2=-2$, $C_2^2=-2$ and $(C_1.C_2)=1$, we can collapse
successively $C_1$ and $C_2$ to get $\P^2$.
The surface $\x[t_0,c_0]$ is obtained by blowing-up 
the ruled surface $Z[t_0,c_0]$ 8 times.
the first blowing-up is one of the blowing-ups described above
according as $c_0 \neq 0$ or $c_0 =0$.
Hence we can blow-down the surface $\x[t_0,c_0]$ to $\P^2$.
\begin{Corollary}
The Picard lattice $\Pic \x[t_0,c_0]$ is isomorphic to the lattice 
$\bigoplus_{i=1}^{10} {\bf Z} e_i $ of rank 10 such that 
\begin{eqnarray*}
(e_i,e_j)=
\left\{ 
\begin{array}{ll}
\delta_{ij} \hspace{1cm} \mbox{if $i=1$},
\smallskip \\
-\delta_{ij} \hspace{1cm} \mbox{otherwise}.
\end{array}
\right.
\end{eqnarray*}
\end{Corollary}

\Proof This is an immediate consequence of Lemma (2.1).
\par
\vspace{1pc}
We have seen in \S1 that on the rational surface $\x[t_0,c_0]$
there exists the divisor $\D[t_0,c_0]$ of which irreducible
components $\D_i[t_0,c_0]$, $0 \leq i \leq 7$ are isomorphic to $\P^1$
with self-intersection number $-2$ and their dual graph is $E_7^{(1)}$
 (cf. Fig.(1.4)).
The linear equivalence classes $[D_i]$'s are linearly independent in the Picard lattice $L=\Pic \x[t_0,c_0]$.
So the Picard lattice $L$ contains a  free ${\bf Z}$-module of rank 8 spanned
by the class $[D_i]$ of $D_i$ in $\Pic \x[t_0,c_0]$ $(0 \leq i \leq 7)$,
of which the bilinear form is given by minus of the Cartan matrix of the affine root system of type $E_7^{(1)}$.
We denote by abuse notation, this subgroup of $L$ with the induced bilinear
form by $L(E_7^{(1)})$.
For a subset $M$ of $L$ we set
\begin{eqnarray*}
M^{\perp}=
\left\{
x\in L\, | \,  (x,y)=0 \hspace{0.5cm} \mbox{for every $x \in M$}
\right\}.
\end{eqnarray*}
Our aim is to prove Theorem (2.24) and (2.31).
We introduced the curves $D_i$, $0 \leq i \leq 7$ on the algebraic surface 
$\x[t_0,c_0]$.
To prove Theorem (2.24), we need some more curves on $\x[t_0,c_0]$.
On the ruled surface $\pi: Z_0[t_0,c_0] \rightarrow \P^1$, 
we have the following curves.
All the irreducible components of the curves are isomorphic to $\P^1$.
{
\addtocounter{equation}{1}
\setcounter{meq}{\value{equation}}
\setcounter{equation}{0}
\renewcommand{\theequation}{\thenewsection.\themeq.\arabic{equation}}
\begin{eqnarray}
\quad  C_1^0[t_0,c_0]:&=&
\overline{\left\{ 
(y_3,z_3) \in W\, | \,  y_3=0
\right\}}
\subset Z_0[t_0,c_0]\\
\smallskip \nonumber\\
&=&
\left\{
(y_3,z_3)\in W_3\, | \,  y_3=0
\right\}
\cup
\left\{
(y_4,z_4) \in W_4\, | \,  y_4=0
\right\},\nonumber
\end{eqnarray}
which is a fiber of $\pi$.
Here $\bar{A}$ denotes the Zariski closure of a subset $A$ of $Z[t_0,c_0]$.
\begin{eqnarray} \quad 
C_2^0[t_0,c_0]:=
\overline{\left\{
(y_3,z_3)\in W\, | \,  y_3z_3-c_0=0
\right\}} \subset Z_0[t_0,c_0].
\end{eqnarray}
so that 
\begin{eqnarray}
\quad 
C_2^0[t_0,c_0]&=&
\left\{
(y_3,z_3)\in W\, | \,  y_3z_3-c_0=0
\right\}
\cup
\left\{
(y_1,z_1) \in W_1\, | \,  z_1=0
\right\}\\
\smallskip \nonumber\\
&&\cup
\left\{(y_4,z_4) \in W_4\, | \,  y_4-c_0z_4=0
\right\}\nonumber
\end{eqnarray}
on $Z_0[t_0,c_0]$. The curve $C_2^0[t_0,c_0]$ is isomorphic to $\P^1$
for $c_0 \neq 0$, and the union of 2 curves isomorphic to $\P^1$ for $c_0=0$.
Namely we have 
\begin{eqnarray}
\quad C_2^0[t_0,c_0]=
\left\{
\begin{array}{ll}
\overline{\left\{(y_1,z_1) \in W_1\, | \,  z_1=0\right\}}, 
\hspace{1cm} &\mbox{if $c_0 \neq 0$},
\\
\smallskip \\
\overline{\left\{
(y_1,z_1) \in W_1\, | \,  z_1=0
\right\}}\cup C_1^0[t_0,c_0],&\mbox{if $c_0 = 0$},
\end{array}
\right.
\end{eqnarray}
on the surface $Z[t_0,c_0]$.
\begin{eqnarray}
\quad C_4^0[t_0,c_0]:=\overline{\left\{
(y_1,z_1)\in W_1\, | \,  2y_1^2+z_1+t_0=0
\right\}} \subset Z_0[t_0,c_0],
\end{eqnarray}
\begin{eqnarray}
\quad C_4^7[t_0,c_0]:&=&\mbox{the proper transform of $C_4^0[t_0,c_0]$ under the 
birational} \\
\smallskip \nonumber\\
&&\mbox{morphism } Z_7[t_0,c_0] \rightarrow Z_0[t_0,c_0].\nonumber
\end{eqnarray}
\begin{eqnarray}
C_5^0[t_0,c_0]:=\overline{\left\{
(y_1,z_1)\in W_1\, | \,  y_1z_1-c_0=0
\right\}} \subset Z[t_0,c_0].
\end{eqnarray}
So
\begin{eqnarray*}
C_5^0[t_0,c_0]&=&
\left\{
(y_1,z_1) \in W_1\, | \,  y_1z_1-c_0=0
\right\}
\cup
\left\{
(y_3,z_3) \in W_3\, | \,  z_3=0
\right\}\\
\smallskip \nonumber\\
&&\cup \left\{(y_2,z_2) \in W_2\, | \,  y_2-c_0z_2=0 \right\}, \nonumber
\end{eqnarray*}
which is isomorphic to $\P^1$ if $c_0 \neq 0$ and, which is the union
of 2 curves isomorphic to $\P^1$ if $c_0=0$.
\begin{eqnarray*}
C_5^0=\overline{\left\{(y_1,z_1) \in W_1\, | \,  y_1=0 \right\}}
 \cup
\overline{\left\{(y_1,z_1) \in W_1\, | \,  z_1=0 \right\}}.
\end{eqnarray*}
\begin{eqnarray}
C_6^0[t_0,c_0]:=
\overline{\left\{
(y_1,z_1)\in W_1\, | \,  2y_1^3+t_0y_1+y_1z_1+c_0+1=0
\right\}}.
\end{eqnarray}
}
\setcounter{equation}{\value{meq}}
$C_6^0[t_0,-1]$ is irreducible if $c_0 \neq -1$ and if $c_0=-1$, it has 
2 irreducible components:
\begin{eqnarray*}
C_6^0[t_0,-1]=
\overline{\left\{ 
(y_1,z_1) \in W_1\, | \,  y_1=0
\right\}}
\cup
C_4^0[t_0, -1].
\end{eqnarray*}
We denote by $C_3[t_0,c_0]$ the exceptional divisor $E_8$ of the last
or the eighth blow-up $\x[t_0,c_0]=Z_8[t_0,c_0] \rightarrow Z_7[t_0,c_0]$
in the construction of $\x[t_0,c_0]$.
We denote the proper transform of $C_j^{(0)}[t_0,c_0]$ by
$C_j[t_0,c_0]$ for $1 \leq j \leq 6$ with $j \neq 3,4$ under the map 
$\x[c_0,t_0] \rightarrow Z_0[t_0,c_0]$. 
$C_4[t_0,c_0]$ is the total transform of $C_4^{(7)}[t_0,c_0]$ under
the birational morphism $\x[t_0,c_0]= Z_8[t_0,c_0] \rightarrow Z_7[t_0,c_0]$.
To simplify
 the notation, the curves $C_j^0[t_0,c_0]$, $C_j[t_0,c_0]$ are respectively 
denoted simply by $C_j^0$, $C_j$ if there is no danger of confusion.
The surface $Z[t_0,c_0] \rightarrow \P^1$ is a ruled surface so  that 
$\Pic Z[t_0,c_0]={\bf Z}C_1^0 \oplus {\bf Z}S $.
We have $\left(C_1^0 \right)^2=-1$, $S^2=2$, $\left( C_1^0.S \right)=1$.
We recall that $S$ is section of the ruled surface $Z[t_0,c_0]$ introduced
in \S1.
\begin{Lemma}
The curve $C_6^0 \cap W_4$ on $W_4$ is defined by
\begin{eqnarray*}
2z_4-y_4^4+t y_4^2 z_4+(2c_0+1) y_4^3z_4=0.
\end{eqnarray*}
This is the curve appeared in Sublemma (3.6) of \cite{Utr}.
\end{Lemma}

\Proof This follows from the definition of the curve $C_6^0$ 
and the coordinate transformations among $W_i$ $(1 \leq i \leq 4)$.
\begin{Lemma}
We have $C_2^0 \sim S-C_1^0$ on $Z[t_0,c_0]$.
\end{Lemma}

\Proof Since $\Pic Z[t_0,c_0] \simeq {\bf Z} C_1^0+{\bf Z} S$,
we can find integers $a$, $b$ such that $C_2^0 \sim a S+b C_1^0$. 
Then 
$1=\left( S.C_2^0 \right)= \left( S.a S+b C_1^0\right)
=aS^2+b(S.C_1^0)=2a+b$, $1=(C_1^0.C_2^0)=(C_1^0.a S+b C_1^0)=a(C_1^0.S)=1$.
Hence $a=1$, $b=-1$.
\begin{Lemma}
Assume $c_0 \neq 0$ so that $C_2^0$ is isomorphic to $\P^1$.
Let $C_1^1$ and $C_2^1$ be the proper transform of $C_1^0$ and $C_2^0$
respectively under the birational morphism
$\pi_1: Z_1[t_0,c_0] \rightarrow Z_0[t_0,c_0]$.
\begin{description}
\item[(i)] $C_1^1$, $C_2^1$ are isomorphic to $\P^1$.
\item[(ii)] $\left(C_1^1.C_2^1 \right)=0$.
\item[(iii)] ${C_1^1}^2={C_2^1}^2=-1$.
\item[(iv)] $\left(C_1^1.E_1\right)=\left( C_2^1. E_1\right)=1$.
\item[(v)] $C_1^1 \cap E_1$, $C_2^1 \cap E_1$ and the center $a_2 \in E_1$
are three distinct points on $E_1$.
\end{description}
\end{Lemma}

\Proof {\bf (i)} follows the fact that the curves $C_1^0$, $C_1^0$ are 
isomorphic
to $\P^1$. The center $a_1$ of the first blowing-up is the point 
$(y_4,z_4)=(0,0)$ on $W_4$. So that $C_1^0$, $C_2^0$ pass through $a_1$.
We have $\left(C_1^0\right)^2=0$ and $\left( C_2^0\right)^2=0$ by Lemma (2.5).
So that $\left(C_1^1\right)^2=\left(C_2^1\right)^2=-1$.
Moreover $C_1^0$ and $C_2^0$ intersect at $a_1$ transversely 
$\left(C_1^1.C_2^1)
\right)=0$.
Easy calculation shows that $C_1^1$ and $C_1^2$ do not pass through 
$a_2 \in E_1 \cap C_6^1$.
Here $C_6^1$ is the proper transform of $C_6^0$ under the birational
map $\pi_1:Z_1[t_0,c_0] \rightarrow Z_0[t_0,c_0]$
(cf. Lemma (2.4) and Sublemma (3.6) of \cite{Utr}).
\begin{Corollary}
If $c_0 \neq 0$, then the configuration of the curves $D_i$ $(0 \leq i \leq 7)$
and $C_1$, $C_2$ is as in Fig.(2.7.1).
We have $(C_1)^2=(C_2)^2=-1$.
\begin{figure}[h]
\unitlength 1.5mm
\begin{center}
\begin{picture}(70,30)

\put(6,2){$D_1$}
\put(0,15){$D_2$}
\put(16,10){$D_3$}
\put(12,25){$D_4$}
\put(30,11){$D_0$}
\put(44,11){$D_5$}
\put(62,15){$D_6$}
\put(50,2){$D_7$}

\put(8,5){\line(0,1){14}}
\put(6,16){\line(1,0){14}}
\put(18,14){\line(0,1){14}}
\put(16,26){\line(1,0){32}}
\put(32,14){\line(0,1){14}}
\put(46,14){\line(0,1){14}}
\put(44,16){\line(1,0){14}}
\put(52,5){\line(0,1){14}}

\put(0,13){\line(1,0){14}}
\put(0,10){\line(1,0){14}}

\put(-5,12){$C_2$}
\put(-5,9){$C_1$}
\end{picture}\\
{\rm Fig.(2.7.1)}
\end{center}
\end{figure}
\end{Corollary}

\Proof In fact, the centers $a_i$ for $2 \leq i \leq 8$ are
infinitely near points of $a_2$.
So $C_1$, $C_2$ are total transform of $C_1^1$ and $C_2^1$
respectively under the morphism $Z_8[t_0,c_0] \rightarrow Z_2[t_0,c_0]$,
which is an isomorphism on a neighborhood of $C_1^1 \cup C_2^1$.

Similarly we can prove the following 
\begin{Lemma}
If $c_0=0$, then $C_2=C_1 \cup C_2'$ such that $C_1$,$C_2'$ are isomorphic to
$\P^1$ and we have $\left( C_1.C_2'\right)=1$, $\left( C_1 \right)^2=-1$,
$\left( C_2' \right)^2=-2$.
The configuration of $C_1$, $C_2'$ is as in Fig.(2.8.1).
\begin{figure}[h]
\unitlength 1.5mm
\begin{center}
\begin{picture}(70,30)

\put(6,2){$D_1$}
\put(2,15){$D_2$}
\put(16,10){$D_3$}
\put(12,25){$D_4$}
\put(30,11){$D_0$}
\put(44,11){$D_5$}
\put(62,15){$D_6$}
\put(50,2){$D_7$}

\put(8,5){\line(0,1){14}}
\put(6,16){\line(1,0){14}}
\put(18,14){\line(0,1){14}}
\put(16,26){\line(1,0){32}}
\put(32,14){\line(0,1){14}}
\put(46,14){\line(0,1){14}}
\put(44,16){\line(1,0){14}}
\put(52,5){\line(0,1){14}}

\put(0,9){\line(1,0){11}}
\put(2,3){\line(0,1){9}}

\put(11,8){$C_1$}
\put(0,0){$C_2'$}
\end{picture}\\
{\rm Fig.(2.8.1)}
\end{center}
\end{figure}

\end{Lemma}

\Proof If $c_0=0$, then a point $(y_1,z_1) \in W_1$ and a point 
$(y_3,z_3) \in W_3$ are identified if $y_1y_3=1$ and $z_3=-y_1z_1$
so that $Z[t_0,0] \simeq F_2$ and we have 2 disjoint sections of 
the ruled surface $Z[t_0,0]$. Namely
\begin{eqnarray*}
S=\left\{
(y_2,z_2) \in W_2\, | \,  z_2=0
\right\}
\cup
\left\{
(y_4,z_4) \in W_4\, | \,  z_4=0
\right\}
\end{eqnarray*}
and
\begin{eqnarray*}
S_0=\left\{
(y_1,z_1) \in W_1\, | \,  z_1=0
\right\}
\cup
\left\{
(y_3,z_3) \in W_3\, | \,  z_3=0
\right\}.
\end{eqnarray*}
Then $S^2=2$, $S_0^2=-2$.
By (2.3.4), we have $C_2^0=S_0 \cup C_1^0$.
$C'_2$ is the total transform of $S_0$ that coincides with its proper 
transform 
because $S_0$ does not pass through the center $a_1$ of the first blowing-up
and consequently any centers $a_j$, $2 \leq j \leq 8$.
In particular $C_2' \subset \x \backslash \left(\bigcup_{i=0}^7 D_i \right)$
(cf. Proposition (4.2) below).

Let us now study the curve $C_4$.
\begin{Lemma}
The curve $C_4^0 \cup W_4$ on the affine plane $W_4 \subset Z_0[t_0, c_0]$
is defined by 
\end{Lemma}
\begin{eqnarray*}
2z_4+(c_0z_4-y_4)y_4^3+ty_4^2z_4=0.
\end{eqnarray*}

\Proof Since 
$C_4^0=
\overline{\left\{
(y_1,z_1) \in W_1 \, | \,  2y_1^2+z_1+t=0
\right\}}$,
this follows from the coordinate transformations among the $W_i$'s
$1 \leq i \leq 4$ given \S1.
\begin{Lemma}
The curve $C_4^0$ has the following properties.
\begin{description}
\item[(i)] The proper transform $C_4^i$ of $C_4^0$ on $Z_i[t_0,c_0]$ 
by the birational morphism $Z_i[t_0,c_0] \rightarrow Z_0[t_0,c_0]$
passes through the center $a_{i+1}$ of the 
$(i+1)$-th blow-up for $0 \leq i \leq 6$.
\item[(ii)]If $c_0 \neq -1$, $C_4^7$ does not pass through $a_8$.
\item[(iii)]If $c_0 =-1$, $C_4^7$ pass through $a_8$.
\end{description}
\end{Lemma}

\Proof This is a consequence of Lemma (2.9) and Sublemma (3.6) of \cite{Utr}.
\begin{Corollary}
\begin{description}
\item[(i)]On $Z[t_0,c_0]$, we have $C_4^0 \sim S+2C_1^0$ so that 
$\left(C_4^0\right)^2=6$. 
\item[(ii)]$\left(C_4^7\right)^2=-1$. 
\item[(iii)]$\left(C_4\right)^2=-1$. 
\item[(iv)]$(C_3.C_4)=0$.
\end{description}
\end{Corollary}

\Proof Let $C_4^0 \sim aS+bC_1^0$ with $a,b \in {\bf Z}$, on $Z_0[t_0,c_0]$.
Since for a given $y_1$, the equation $2y_1^2+z_1+t=0$ for $z_1$
has a unique solution,
\begin{eqnarray*}
1=\left(C_1^0\,.\,C_4^0\right)=\left(C_1^0\,.\,aS+bC_1^0 \right)=a.
\end{eqnarray*}
The defining equation for $C_4^0$ on $W_2$ is 
$2y_2^2z_2+1+tz_2=0$.
So $C_4^0 \cap W_2 \cap S=\emptyset$.
The defining equation of $C_4^0$ on $W_4$ is 
\begin{eqnarray*}
2z_4+(c_0z_4-y_4)y_4^3+ty_4^2z_4=0
\end{eqnarray*}
so that on $W_4$, $W_4 \cap C_4^0 \cap S=(0,0) \in W_4$.
To multiplicity of $W_4 \cap C_4^0 \cap S$ is 4.
Hence 
\begin{eqnarray*}
4=\left(S\,.\,C_4^0 \right)=\left(S\,.\,aS+bC_1^0\right)=2a+b
\end{eqnarray*}
so that $b=2$. Consequently $C_4^0 \sim S+2C_1^0$.
\begin{eqnarray*}
\left(
C_4^0
\right)^2=\left(S+2C_1^0\right)^2=S^2+4\left(S\,.\,C_1^0\right)+4\left(C_1^0\right)^2=2+4=6.
\end{eqnarray*}
Therefore $\left( C_4^7\right)^2=-1$ by Lemma (2.10).
(iii) and (iv) follow from (ii) and Lemma (2.10).

\begin{Corollary}
If $c_0 \neq -1$,
{\bf (i)} $C_4$ is an irreducible curve,
{\bf(ii)} we have $\left(C_3\right)^2=\left( C_4\right)^2=-1$ and
{\bf (iii)} the configuration of the curves $D_i$ $(0 \leq i \leq 7)$
and $C_3$, $C_4$ is as Fig.(2.12.1).
\unitlength 1.5mm
\begin{figure}[h]
\begin{center}
\begin{picture}(70,30)

\put(6,2){$D_1$}
\put(0,15){$D_2$}
\put(16,10){$D_3$}
\put(12,25){$D_4$}
\put(30,11){$D_0$}
\put(44,11){$D_5$}
\put(62,15){$D_6$}
\put(52,2){$D_7$}

\put(8,5){\line(0,1){14}}
\put(6,16){\line(1,0){14}}
\put(18,14){\line(0,1){14}}
\put(16,26){\line(1,0){32}}
\put(32,14){\line(0,1){14}}
\put(46,14){\line(0,1){14}}
\put(44,16){\line(1,0){14}}
\put(56,5){\line(0,1){14}}

\put(51,13){\line(1,0){14}}
\put(51,9){\line(1,0){14}}

\put(66,12){$C_3$}
\put(66,8){$C_4$}

\end{picture}\\
{\rm Fig.(2.12.1)}
\end{center} 
\end{figure}

\end{Corollary}

\Proof This is a direct consequence of Lemma (2.10) and Corollary (2.11).

\begin{Corollary}
If $c_0=-1$, then $C_4=C_3 \cup C_4'$, where $C_4'$ is isomorphic to $\P^1$.
We have $\left( C_3\right)^2=-1$, $\left( C_4'\right)^2=-2$, 
$\left( C_3\,.\,C_4'\right)=1$.
We also have $\left( C_3 \right)^2=\left( C_4\right)^2=-1$.
$C_4'$ coincides with a component of $C_6$
\end{Corollary}

\Proof The decomposition of $C_4$ follows from Lemma (2.10).
Since $C_3$ is the exceptional curve of the first kind on $\x[t_0,c_0]$,
$\left(C_3\right)^2=-1$.
We can check $(C_3\,.\,C_4')=1$ by a direct calculation.
By Corollary (2.11), (ii),
\begin{eqnarray*}
-1=\left( C_4 \right)^2=\left(C_3+C_4' \right)^2
=\left( C_3 \right)^2+2\left( C_3\,.\,C_4' \right)+\left( C_4' \right)^2
=-1+2+\left( C_4' \right)^2=1+\left( C_4' \right)^2.
\end{eqnarray*}
So $\left( C_4' \right)^2=-2$.
The last assertion follows from the definition of $C_4'$ and $C_6$.
\begin{Corollary}
If $c_0=-1$, the configuration of the curves $D_i$ $(0\leq i \leq 7)$
and $C_3$, $C_4'$ is as in Fig.(2.14.1).
\begin{figure}[h]
\unitlength 1.5mm
\begin{center}
\begin{picture}(70,30)

\put(6,2){$D_1$}
\put(2,15){$D_2$}
\put(16,10){$D_3$}
\put(12,25){$D_4$}
\put(30,11){$D_0$}
\put(44,11){$D_5$}
\put(62,15){$D_6$}
\put(50,2){$D_7$}

\put(8,5){\line(0,1){14}}
\put(6,16){\line(1,0){14}}
\put(18,14){\line(0,1){14}}
\put(16,26){\line(1,0){32}}
\put(32,14){\line(0,1){14}}
\put(46,14){\line(0,1){14}}
\put(44,16){\line(1,0){14}}
\put(52,5){\line(0,1){14}}

\put(49,9){\line(1,0){11}}
\put(57,3){\line(0,1){9}}

\put(61,8){$C_3$}
\put(55,0){$C_4'$}
\end{picture}\\
{\rm Fig.(2.14.1)}
\end{center}
\end{figure}
\end{Corollary}
\begin{Lemma} We have 
{\bf (i)} $\left( C_1\,.\,C_3\right)=0$, 
{\bf (ii)} $\left( C_1\,.\,C_4\right)=0$,
{\bf (iii)} $\left(C_2\,.\, C_3 \right)=0$ and
{\bf (iv)} $\left(C_2\,.\,C_4 \right)=2$.
\end{Lemma}

\Proof
By Lemma (2.5) and Corollary (2.11)
\begin{eqnarray*}
C_2^0 \sim S-C_1^0, \hspace{1cm} C_4^0 \sim S+2C_1^0
\end{eqnarray*}
so that 
\begin{eqnarray*}
\left(C_2^0\,.\, C_4^0 \right)=\left(S-C_1^0\,.\, S+2C_1^0 \right)
=S^2+(C_1^0\, .\, S)=2+1=3.
\end{eqnarray*}
Both $C_2^0$ and $C_4^0$ pass through the point $(y_4,z_4)=(0,0)$ on $W_4$,
which is not a common point of $C_2^1$ and $C_4^1$.
Here $C_i^1$ is the proper transform of $C_i^0$ for birational morphism
$\pi_1: Z_1[t_0,c_0] \rightarrow Z_0[t_0,c_0]$.
Hence $2=\left(C_2^1\,.\,C_4^1 \right)=\left(C_2\,.\,C_4 \right)$. 
Other intersection
numbers are evident from the geometric definition of the curves.
\par
\vspace{1pc}
As we imagine from Fig.(2.7.1) and (2.12.1), (2.8.1) and (2.14.1).
The pairs $\{C_1,C_2\}$, $\{C_3,C_4\}$ of curves are symmetric with respect to 
$J_{\x}(t_0;c_0,-1-c_0)$ (cf. Proposition (1.10)).
\begin{Proposition}
\begin{description}
\item[(i)]We have 
\begin{eqnarray*}
J_{\x}(t_0;c_0,-1-c_0)\left(C_1[t_0,c_0] \right)=C_3[t_0,-1-c_0],\\
\smallskip \nonumber\\
J_{\x}(t_0;c_0,-1-c_0)\left(C_2[t_0,c_0] \right)=C_4[t_0,-1-c_0],\\
\smallskip \nonumber\\
J_{\x}(t_0;c_0,-1-c_0)\left(C_3[t_0,c_0] \right)=C_1[t_0,-1-c_0],\\
\smallskip \nonumber\\
J_{\x}(t_0;c_0,-1-c_0)\left(C_4[t_0,c_0] \right)=C_2[t_0,-1-c_0].
\end{eqnarray*}
\item[(ii)]Moreover we have 
\begin{eqnarray*}
J_{\x[t_0,c_0]}(c_0,-1-c_0)\left(C_5[t_0,c_0] \right)=C_6[t_0,-1-c_0],\\
\smallskip \nonumber\\
J_{\x[t_0,c_0]}(c_0,-1-c_0)\left(C_6[t_0,c_0] \right)=C_5[t_0,-1-c_0].
\end{eqnarray*}
\end{description}
\end{Proposition}

\Proof These formulas follow from the definition of curves and from
the proof of Lemma (3.5)
in \cite{Utr}, where the local forms $J_{ij}:W_i \rightarrow W_j$
of $J_{\x}[t_0;c_0,-1-c_0] $ are given. 
\par
\vspace{1pc}
We can prove the following proposition in a similar way.
\begin{Proposition}
\begin{description}
\item[(i)] If $c_0 \neq 0$, we have 
\begin{eqnarray*}
I_{\x}[t_0;c_0,-c_0]\left(C_2[t_0,c_0]\right)&=&C_1[t_0,-c_0],\\
\smallskip \nonumber\\
I_{\x}[t_0;c_0,-c_0]\left(C_1[t_0,c_0]\right)&=&C_2[t_0,-c_0].
\end{eqnarray*}
\item[(ii)] We also have 
\begin{eqnarray*}
I_{\x}[t_0;c_0,-c_0]\left( C_3[t_0,c_0]\right)=C_3[t_0, -c_0].
\end{eqnarray*}
\end{description}
\end{Proposition}
\begin{Remark}
If $c_0=0$, $I_{\x}[t_0,c_0]=Id_{\x[t_0,c_0]}$.
If $c_0\neq -1 $, then
$$
C_4'':=I_{\x}[t_0;c_0,-c_0]\left( C_4[t_0,c_0]\right)
$$
is a new $-1$-curve on $I_{\x}[t_0,-c_0]$
with $\left( C_4''[t_0,-c_0]\,.\,C_1[t_0,-c_0]\right)=2$. 
\end{Remark}

We have to study the curves $C_5$ and $C_6$. We have 
the following lemmas.
\begin{Lemma}
We have 
${\bf(i)}\left(C_5\,.\,D_0\right)=1$,
${\bf(ii)}\left(C_5\,.\,C_1\right)=1$,
${\bf(iii)}\left(C_5\,.\,C_2\right)=0$,
${\bf(iv)}\left(C_5\,.\,C_3\right)=0$,
${\bf(v)}\left(C_5\,.\,C_4\right)=3$,
${\bf(vi)}\left(C_5\,.\,D_i\right)=0$ for $2 \leq i \leq 7$,
${\bf(vii)}\left(C_5\,.\,D_1\right)=1$.
\end{Lemma}

\begin{Lemma}
We have 
${\bf(i)}\left(C_6\,.\,D_0\right)=1$,
${\bf(ii)}\left(C_6\,.\,C_1\right)=0$,
${\bf(iii)}\left(C_6\,.\,C_2\right)=3$,
${\bf(iv)}\left(C_6\,.\,C_3\right)=1$,
${\bf(v)}\left(C_6\,.\,C_4\right)=0$,
${\bf(vi)}\left(C_6\,.\,D_i\right)=0$ for $1 \leq i \leq 6$,
${\bf(vii)}\left(C_6\,.\,D_7\right)=1$.
\end{Lemma}

We prove only Lemma (2.20).
Then Lemma (2.19) is proved
by a similar method or once we have Lemma (2.20),
by the isomorphism
$
J_{\x[t_0,c_0]}(c_0,-1-c_0):\x[t_0,c_0] \rightarrow \x[t_0, -1-c_0]
$
in view of Proposition (1.11) and (2.16).

\Proof
The curve $C_6^0 \cap W_4$ on $W_4$ is defined by
\addtocounter{equation}{1}
\setcounter{meq}{\value{equation}}
\setcounter{equation}{0}
{
\renewcommand{\theequation}{\thenewsection.\themeq.\arabic{equation}}
\begin{eqnarray}
2z_4-y_4^4+ty_4^2z_4+(2c_0+1)y_4^3z_4=0
\end{eqnarray}
by Lemma (2.4).
$C_6^0$ and 
\begin{eqnarray}
C_1^0 \cap W_4 =\left\{(y_4,z_4) \in W_4\, | \,  z_4=0 \right\}
\end{eqnarray}
}\setcounter{equation}{\value{meq}}
have a point $(y_4,z_4)=0$ in common.
Since by (2.19.1) and (2.19.2),
\begin{eqnarray*}
z_4=\frac{y_4^4}{2+ty_4^2+(2c_0+1)y_4^3}.
\end{eqnarray*}
The intersection multiplicity of $C_6^0 \cap S$ at $(y_4,z_4)=(0,0)$ is 4.
On the other hand, the curve $C_6^0 \cap W_2$ on $W_2$ is defined by
$2y_2^3z_2+y_2+ty_2z_2+(c_0+1)z_2=0$ by (2.3.8),
which is definition of curve $C_6^0$.
So the intersection multiplicity $C_6^0 \cap S$ at
$(y_2,z_2)=(0,0)$ is 1. Therefore we have 
\addtocounter{equation}{1}
\setcounter{meq}{\value{equation}}
\setcounter{equation}{0}
{
\renewcommand{\theequation}{\thenewsection.\themeq.\arabic{equation}}
\begin{eqnarray}
\left( C_6^0\,.\,S\right)=5.
\end{eqnarray}
It follows from (2.19.1) and (2.3.1)
\begin{eqnarray}
\left( C_6^0\,.\,C_1^0\right)=1.
\end{eqnarray}
}\setcounter{equation}{\value{meq}}
Writing $C_6^0 \sim aS+bC_1^0$, we conclude from (2.20.1)
and (2.20.2), 
\begin{eqnarray}
C_6^0 \sim S+3 C_1^0.
\end{eqnarray}
In particular we have $\left(C_6^0 \right)^2=8$.
Now it follows from Sublemma (3.6) of \cite{Utr}, $\left(C_6\right)^2=0$
because $C_6^0$ passes through all the centers $a_i$ $(1 \leq i \leq 8)$.
Now $\left( C_6^0\,.\,C_2^0 \right)=(S+3C_1^0\,.\,S-C_1^0)
=S^2+2(C_1^0\,.\,S)=4$.
Since the intersection multiplicity of $C_6^0 \cap C_2^0$
at $(y_4,z_4) =0$ is 1, we have $(C_6\,.\,C_2)=3$.
Other multiplicities are calculate in a similar way starting from (2.23).
\begin{Theorem}
{\bf (i)} $L(A_1^{(1)}):=L\left( E_7^{(1)}\right)^{\perp}$ is isomorphic to 
${\bf Z}^{\oplus 2}$ with the inner product given by minus of Cartan matrix
\begin{eqnarray*}
\left[
\begin{array}{cc}
-2&2\\
\smallskip \\
2&-2
\end{array}
\right]
\end{eqnarray*}
of the affine root system of type $A_1^{(1)}$.\\
{\bf(ii)} We have 
$L\left( A_1^{(1)}\right)^{\perp}=L\left( E_7^{(1)}\right)$.
\end{Theorem}

\Proof We denote by $[C]$ the linear equivalence class of a curve,
which we often denote by abuse of notation simply by $C$.
We set 
\begin{eqnarray*}
M=\left< [C_2]-[C_1], [C_4]-[C_3]\right> \subset \Pic \x[t_0,c_0].
\end{eqnarray*}
So $M$ is a subgroup of rank 2 in $\Pic \x[t_0,c_0]$.
We have $\left( C_2-C_1\right)^2=\left( C_2\right)^2 + \left( C_1\right)^2=-2$,
$\left( C_4 -C_3\right)^2=\left( C_4 \right)^2 + \left( C_3\right)^2 =-2 $
by Corollary (2.7), Lemma (2.8), Corollaries (2.12) and (2.14).
Moreover
$\left( C_2-C_1\,.\, C_4-C_3\right)=\left(C_2\,.\,C_4\right)
-\left( C_1\,.\,C_4\right)
-\left(C_2\,.\,C_3\right)+\left( C_1\,.\,C_3\right)=2$
by Lemma (2.4).
Therefore the inner product on $M$ is given by minus of the Cartan matrix of the affine root system of type $A_1^{(1)}$.
It follows from Lemma (2.6), (2.8), Corollaries (2.12), (2.14) 
$M \subset L\left(E_7^{(1)} \right)^{\perp}$.
We have to show $M=L\left(E_7^{(1)}\right)^{\perp}$.
Since $M$ and $L\left(E_7^{(1)}\right)^{\perp}$ are of rank 2, we have 
${\mathbb Q} \otimes _{\bf Z}M
={\mathbb Q} \otimes _{\bf Z} L\left(E_7^{(1)}\right)^{\perp}$.
We work in ${\mathbb Q} \otimes _{\bf Z} \Pic \x[t_0,c_0]$.
Let
\begin{eqnarray*}
a\left(C_2-C_1\right)+b\left( C_4-C_3\right) \in {\mathbb Q} \otimes _{\bf Z} M
\subset {\mathbb Q} \otimes _{\bf Z} \Pic \x[t_0,c_0]
\hspace{1cm}\mbox{$(a,b \in \mathbb Q)$}.
\end{eqnarray*}
To prove (i), we have to show that if
\begin{eqnarray}
a\left(C_2-C_1\right)+b\left( C_4-C_3\right)  \in \Pic \x[t_0,c_0],
\end{eqnarray}
then $a,b \in {\bf Z}$.
In fact, by (2.23)
\addtocounter{equation}{1}
\setcounter{meq}{\value{equation}}
\setcounter{equation}{0}
{
\renewcommand{\theequation}{\thenewsection.\themeq.\arabic{equation}}
\begin{eqnarray}
\left(C_2\,.\,a\left(C_2-C_1\right)+b\left(C_4-C_3\right) \right)
=-a+2b\in {\bf Z},\\
\smallskip \nonumber\\
\left(C_4\,.\,a\left(C_2-C_1\right)+b\left(C_4-C_3\right) \right)
=2a-b\in {\bf Z},\\
\smallskip \nonumber\\
\left(C_5\,.\,a\left(C_2-C_1\right)+b\left(C_4-C_3\right) \right)
=-a+3b\in {\bf Z},\\
\smallskip \nonumber\\
\left(C_6\,.\,a\left(C_2-C_1\right)+b\left(C_4-C_3\right) \right)
=3a-b\in {\bf Z}.
\end{eqnarray} 
}
\setcounter{equation}{\value{meq}}
It follows from (2.26.1) and (2.26.3), $b \in {\bf Z}$.
Similarly by (2.26.2) and (2.26.4) we conclude $a \in {\bf Z}$.
Now we show $L\left(A_1^{(1)} \right)^{\perp}=L\left( E_7^{(1)}\right)$.
Since $L\left( E_7^{(1)}\right) \subset L\left( A_1^{(1)}\right)^{\perp}$,
we may argue as above.
Namely $F=\sum_{i=0}^{7} a_i D_i 
\in {\mathbb Q} \otimes _{\bf Z} L\left( E_7^{(1)}\right)$ with 
$a_1 \in {\mathbb Q}$,
we have to show that
if $F \in \Pic \x[t_0,c_0]$, then $a_i \in {\bf Z}$.
Since $(C_1 \,.\, F) \in {\bf Z}$ 
which is equal to $a_i$ by Corollary (2.7) and 
Lemma (2.8), so $a_1 \in {\bf Z}$.
Similarly using $C_3$, we conclude $a_7 \in {\bf Z}$.
Now $\left( C_5 \,.\, F\right) \in {\bf Z}$ which is equal to $a_1+a_0$ by
Lemma (2.19) so that $a_1+a_0 \in {\bf Z}$ and consequently $a_0 \in {\bf Z}$.
Now since 
\begin{eqnarray*}
\left(D_1\,.\,F \right)=-2a_1+a_2 \in {\bf Z},\\
\smallskip \nonumber\\
\left(D_7\,.\,F \right)=-2a_7+a_6 \in {\bf Z},
\end{eqnarray*}  
we conclude $a_2,a_6 \in {\bf Z}$.
We have further
\begin{eqnarray*}
\left(D_2\,.\,F \right)=a_1-2a_2+a_3 \in {\bf Z},\\
\smallskip \nonumber\\
\left(D_6\,.\,F \right)=a_7-2a_6+a_5 \in {\bf Z},
\end{eqnarray*}
so that $a_3,a_5 \in {\bf Z}$.
Finally 
\begin{eqnarray*}
\left(D_3\,.\,F \right)=a_2-2a_3+a_4 \in {\bf Z}
\end{eqnarray*}
so that $a_4 \in {\bf Z}$.

Writing $\x[t_0,c_0]$ by $X$,
we have an exact sequence of homology
\begin{eqnarray}
\rightarrow H_3(X;D,{\bf Z})
\rightarrow H_2(X \backslash D,{\bf Z})
\rightarrow H_2(X,{\bf Z})
\rightarrow H_2(X;D,{\bf Z})
\rightarrow .
\end{eqnarray}
By the Poincar\'e duality, we have 
\begin{eqnarray}
H_3(X;D,{\bf Z})=H^1(D, {\bf Z})=0
\end{eqnarray}
and 
\begin{eqnarray}
H_2(X;D,{\bf Z})=H^2(D,{\bf Z}) \simeq \bigoplus_{i=0}^{7} {\bf Z} D_i.
\end{eqnarray}
Since $X$ is a non-singular, projective rational surface 
$\Pic X \simeq H_2(X, {\bf Z})$.
So it follows from the exact sequence (2.25)
and (2.28), (2.29) that 
\begin{eqnarray}
L\left( E_7^{(1)}\right)^{\perp} = H_2(X \backslash D, {\bf Z}).
\end{eqnarray}
Let $i:\C(t,c,) \rightarrow \C(t,c)$ be the $\C(t)$-automorphism
of the field $\C(t,c)$ sending $c$ to $-c$.
Similarly, let $j:\C(t,c) \rightarrow \C(t,c)$ be the $\C(t)$-automorphism
of the field $\C(t,c)$ such that $j(c)=-1-c$. So the 
subgroup $<i,j>$ of the automorphisms group of the field $\C(t,c)$
generated by $i$ and $j$ is isomorphic to the extended affine Weyl group
$\widehat{G}$.
The group $W_a$ operates on $\Spec \C (t,c)$.
The group $\widehat{G}$ also operates on 
$\x_{\C(t,c)}:=\x \otimes_{\C[t,c]} \C(t,c)$
in such a way that the morphism $\x_{\C(t,c)} \rightarrow \Spec \C(t,c)$
is $\widehat{G}$-equivariant.
A $-1$-curve $C$ on a surface is a curve isomorphic to $\P^1$ with $C^2=-1$.
\begin{Theorem}
For every $(t_0,c_0) \in \C^2$, there are infinitely many $-1$-curves on 
$\x[t_0,c_0]$.
\end{Theorem}

\Proof
First we assume that $c_0$ is not an integer.
Then we have an isomorphism
\begin{eqnarray*}
T_+[t_0\,;\,c_0-1\,,\,c_0]=I_{\x}[t_0\,;\,-c_0\,,\,c_0] 
\circ J_{\x}[t_0\,;\,c_0-1\,,\,-c_0]
:\x[t_0\,,\,c_0-1] \rightarrow \x[t_0\,,\,c_0].
\end{eqnarray*}
So in particular the isomorphism 
\begin{eqnarray*}
T_{+}[t_0\,;\,c_0-1\,,\,c_0]:X[t_0\,,\, c_0-1] \rightarrow \x[t_0\,,\, c_0]
\end{eqnarray*}
maps the $-1$-curve $C_3[t_0\,,\,c_0-1]$ on $\x[t_0\,,\,c_0-1]$
to a $-1$-curve 
$$
\Gamma_1[t_0\,,\,c_0]:=T_{+}[t_0\,;\,c_0-1\,,\,c_0](C_3[t_0\,,\,c_0-1])
$$
on $\x[t_0 \,,\, c_0]$.
For a positive integer $n$, the iteration
\begin{eqnarray*}
\lefteqn{T_{+}[t_1\,;\,c_0-n\,,\,c_0]:=
T_{+}[t_0\,;\,c_0-1\,,\,c_0]
\circ 
T_{+}[t_0\,;\,c_0-2\,,\,c_0-1]}
\\ \smallskip \\
&&\circ
\cdots
\circ
T_{+}[t_0\,;\,c_0-n+1\,,\,c_0-n+2]
\circ
T_{+}[t_0\,;\,c_0-n\,,\,c_0-n+1]\\
\smallskip \\
&&:\x[t_0\,,\,c_0-n] \rightarrow \x[t_0\,,\,c_0-n+1] 
\rightarrow \cdots \rightarrow
\x[t_0\,,\,c_0-1] \rightarrow \x[t_0\,,\,c_0]
\end{eqnarray*}
of isomorphisms maps the $-1$-curve $C_3[t_0 \,,\, c_0-n]$ to a $-1$-curve
$\Gamma_n[t_0\,,\,c_0]:=T_{+}[t_0 \,;\, c_0-n,c_0](C_3[t_0 \,,\,c_0-n])$.
Hence we have $-1$-curves $\Gamma_n[t_0\,,\,c_0]$ $(n=1,2,\ldots)$,
on $\x[t_0 \,,\, c_0]$.
We have to show that the $-1$-curves $\Gamma_n[t_0\,,\,c_0]$ are 
distinct curves on $\x[t_0\,,\,c_0]$.
To this end, it is sufficient to show that the linear equivalence classes
$[\Gamma_n] \in \Pic \x[t_0\,,\,c_0]$ are distinct.
It follows from the construction of $\x[t_0,c_0]$ that
\begin{eqnarray*}
\Pic \x[t_0\,,\, c_0] = 
\bigoplus_{i=0}^{7} {\mathbb Z} D_i 
\bigoplus {\mathbb Z} C_1
\bigoplus {\mathbb Z} C_3.
\end{eqnarray*}
By Proposition (2.16) and (2.17), we have 
\begin{eqnarray}
T_{+}[t_0 \,;\, c_0-1 \,,\, c_0]
\left(C_1[t_0\,,\,c_0-1]\right)=C_3[t_0 \,,\, c_0],
\nonumber\\
\smallskip \nonumber\\
T_{+}[t_0 \,;\, c_0-1 \,,\, c_0]
\left(C_3[t_0\,,\,c_0-1]\right)=C_2[t_0 \,,\, c_0].
\end{eqnarray}
Let us express the linear equivalence class of $C_2$ as a linear combination
of $C_1$, $C_3$ and the $D_i$'s.
Namely we set 
\begin{eqnarray*}
C_2 \sim a_1 C_1+a_3 C_3+ \sum_{i=0}^{7} b_i D_i
\end{eqnarray*}
and determine the integers $a_1$, $a_3$ and $b_j$'s.
We have linear equations
\begin{eqnarray*}
&&0=\left( C_1 \,,\, C_2 \right)=-a_1+b_1,\\
\smallskip \\
&&1=\left( D_1 \,,\, C_2 \right)=a_1-2 b_1+b_2,\\
\smallskip \\
&&0=\left( D_2 \,,\, C_2 \right)=b_1-2 b_2+b_3,\\
\smallskip \\
&&0=\left( D_3 \,,\, C_2 \right)=b_2-2 b_3+b_4,\\
\smallskip \\
&&0=\left( D_4 \,,\, C_2 \right)=b_3-2 b_4+b_5+b_0,\\
\smallskip \\
&&0=\left( D_0 \,,\, C_2 \right)=-2 b_0+b_4,\\
\smallskip \\
&&0=\left( D_5 \,,\, C_2 \right)=b_4-2 b_5+b_6,\\
\smallskip \\
&&0=\left( D_6 \,,\, C_2 \right)=b_5-2 b_6+b_7,\\
\smallskip \\
&&0=\left( D_7 \,,\, C_2 \right)=-2 b_7+b_6+a_3,\\
\smallskip \\
&&0=\left( C_3 \,,\, C_2 \right)=-b_7-a_3.
\end{eqnarray*}
We solve this system of linear equations to get
\begin{eqnarray*}
\Gamma_1=C_2 \sim -C_1+2 C_3 -D_1+D_3+2 D_4 +2 D_5 +2 D_6 +2 D_7 + D_0.
\end{eqnarray*}
In particular, we have 
\begin{eqnarray}
&&\Gamma_1=T_+[t_0 \,;\, c_0-1 \,,\, c_0](C_3) \equiv  -C_1+2 C_3 \\
\smallskip \nonumber \\
&&T_+[t_0 \,;\, c_0-1 \,,\, c_0](C_1) \equiv  C_3 
\end{eqnarray}
modulo the subgroup 
$L\left( E_7^{(1)} \right)=\left< D_0 \,,\,D_1\,,\, \ldots \,,\, D_7 \right>$
in $\Pic X[t_0 \,,\, c_0]$.
Since $L\left( E_7^{(1)} \right)$ is invariant under $I_{\x}$ and $J_{\x}$
by Proposition (1.10) and hence by $T_{+}$, it follows from (2.34) and (2.35),
\begin{eqnarray*}
\Gamma_2&=&T_{+}[t_0\,;\,c_0-1\,,\,c_0] \circ T_{+}[t_0 \,;\, c_0-2 \,,\,c_0-1]
\left( C_3 [t_0\,,\, c_0-2]\right)\\
\smallskip \\
&\equiv& T_{+}[t_0\,;\,c_0-1\,,\,c_0] 
\left( -C_1 [t_0\,,\, c_0-1]+2C_3[t_0 \,;\,c_0-1]\right)\\
\smallskip \\
&\equiv& -C_3[t_0 \,,\, c_0]+
2 \left( -C_1 [t_0\,,\, c_0]+2C_3[t_0 \,;\,c_0]\right)\\
\smallskip \\
&\equiv& -2 C_1 [t_0\,,\, c_0]+3C_3[t_0 \,;\,c_0]
\end{eqnarray*}
modulo $\Gamma\left( E_7^{(1)} \right)$.
So $\Gamma_1$ and $\Gamma_2$ are distinct curves on $X[t_0\,,\,c_0]$. 
Similarly we have 
\begin{eqnarray*}
\Gamma_3 &\equiv& -6 C_1+ 10 C_3\\
\smallskip \\
\Gamma_4 &\equiv& -20 C_1+ 34 C_3\\
\smallskip \\
\Gamma_5 &\equiv& -34 C_1+ 116C_3\\
\smallskip \\
&&\vdots
\end{eqnarray*} 
so that $\Gamma_1$, $\Gamma_2$, $\cdots$ are distinct $-1$-curves on 
$\x[t_0\,,\,c_0]$. 
If $c_0$ is an integer, since the surface $\x[t_0 \,,\, c_0]$
is isomorphic to $\x[t_0,-1]$ as we proved in \cite{Utr},
we may assume $c_0=-1$. Then we can apply the above argument.
\section{Vanishing cycles on $\x[t_0\,,\,c_0]$.}
\setcounter{equation}{0}
\par
If $c_0 =1$, then we have the $-2$-curve $C_2'[t_0\,,\,0]$
on $\x[t_0,0]$ that is the limit of the cycle
$C_2[t_0,c]-C_1[t_0,c]$ on $\x[t_0,c]$ when $c \rightarrow 0$.
We can collapse the $-2$-curve $C_2'[t_0,0]$ to a normal singular point to get
a normal algebraic surface $\x''[t_0]$ with an ordinary double point.
This suggests that the cycle $C_2[t_0,c]-C_1[t_0,c]$ on $\x[t_0,c]$ is a
vanishing cycle.
\begin{Proposition}
There exists a flat family 
$\varphi':\x'[t,c^2] \rightarrow \Spec \C[t,c^2]$ of projective algebraic
surfaces with the following properties.\\
{\bf (i)} The variety $\x'[t,c^2]$ is non-singular.\\
{\bf (ii)} We have an isomorphisms 
\begin{eqnarray*}
\x[t_0,c_0] \simeq \x'[t_0,c_0^2] \hspace{0.5cm} \mbox{ if $c_0 \neq 0$,}
\end{eqnarray*}
and 
\begin{eqnarray*}
\x''[t_0] \sim \x'[t_0,0].
\end{eqnarray*}
Here $\x'[t_0,c_0^2]$ denotes the fiber ${\varphi'}^{-1}((t_0,c_0^2))$
over a point $(t_0,c_0^2) \in \Spec \C[t,c^2]$.\\
{\bf (iii)} 
The morphism $\varphi:\x'[t,c^2] \rightarrow \Spec \C[t,c^2]$ is
smooth except for the points on $\x'[t_0,0]$
that correspond to the rational double point of $\x''[t_0]$, $t_0 \in \C$.\\
{\bf (iv)}
We have a $\C[t,c,1/c]$-isomorphism
\begin{eqnarray*}
\x[t,c] \otimes_{\C[t,c]} \C\left[t,c,\frac{1}{c}\right] \simeq 
\x'[t,c^2] \otimes_{\C[t,c^2]} \C\left[t,c, \frac{1}{c}\right].
\end{eqnarray*}
\end{Proposition}

\Proof
We consider a $\C[t,c]$-morphism
\begin{eqnarray*}
f_1&:&W_1 \times \Spec \C[t,c] = W_1 \times \C^2 \rightarrow \P^3 \times \C^2\\
\smallskip \nonumber \\
&&(y_1,t_1;t,c) \mapsto (1,y_1(c-y_1z_1), 2y_1z_1-c,z_1;t,c).
\end{eqnarray*}
The morphism $f_1$ extends to a rational map
\begin{eqnarray*}
f:Z_0[t,c] \rightarrow \P^3 \times \C^2,
\end{eqnarray*}
which turns out as we can check it easily, to be a morphism.
We recall that $Z_0[t,c]$ is the starting family of ruled surfaces in the 
construction of $\x[t,c]$.
For example, we have on $W_3 \times \Spec \C[t,c]$,
\begin{eqnarray*}
f_3&:&W_3 \times \C^2 \rightarrow \P^3 \times \C^2\\
\smallskip \nonumber \\
&&(y_3,z_3;t,c) \mapsto (1,z_3,c-2y_3z_3, y_3(c-y_3z_3);t,c).
\end{eqnarray*}
Moreover $f_3$ factors through a family of quadratic surface 
\begin{eqnarray*}
Q[t,c^2]:4x_1x_3-x_2^2+c^2x_0^2=0
\end{eqnarray*}
in $\P^3 \times \C^2$
so that the surface $Q[t,c^2]$ is defined over $\C[t,c^2]$.
If $c=0$, $f_3$ collapses the curve $C_2'[t_0,0]$ for every $t_0$.
We can check that the morphism $f:Z_0[t,c] \rightarrow Q[t,c^2]$ is an
isomorphism outside the curve $C_2'[t,0]$.
We can also show by calculation that we can blow up $Z_0[t,c]$ and $Q[t,c^2]$
so that we have a morphism 
\begin{eqnarray*}
\widetilde{f}:\x[t,c] \rightarrow \widetilde{Q}[t,c^2]=:\x'[t,c^2]
\end{eqnarray*} 
(cf. Proposition (1.11)).
Then $\widetilde{f}$ has the required properties.
We proved in \cite{Utr} that we have an isomorphism
\begin{eqnarray*}
I_\x[t;c,-c]:\x[t,c] \otimes_{\C[t,c]} \C[t,c,1/c] \rightarrow 
\x[t,c] \otimes_{\C[t,c]} \C[t,c,1/c]
\end{eqnarray*}
covering $\Spec \C[t,c]=\C^2 \rightarrow \Spec[t,c]=\C^2$, 
$(t,c) \mapsto (t,-c)$.
The quotient variety
\begin{eqnarray*}
&\x[t,c] \otimes_{\C[t,c]} \C[t,c,1/c] \slash \left<I_\x[t;c,-c] \right>&\\
&\downarrow& \nonumber \\
&\Spec \C[t,c,1/c]^{<I>}=\Spec \C[t,c^2,1/c^2] &
\end{eqnarray*}
is isomorphic to
\begin{eqnarray*}
&\x'[t,c^2] \otimes_{\C[t,c^2]} \C[t,c^2,1/c^2]&\\
&\downarrow& \nonumber\\
&\Spec \C[t,c^2,1/c^2] .&
\end{eqnarray*}
This shows that the involution $I_\x[t;c,-c]$ arises form the rational double
point of type $A_1$.
Namely if we fix $t=t_0 \in \C$,
the birational involution
\begin{eqnarray}
I:\x[t_0,c] \rightarrow \x[t_0,c]
\end{eqnarray}
comes from the automorphism of the singular quadratic surface
\begin{eqnarray*}
Q=\left\{ 
(x_0,x_1,x_2;c) \in \P^2 \times \C
\, | \,  
4x_1x_3-x_2^2+c^2x_0=0
\right\}
\end{eqnarray*}
sending $(x_0,x_1,x_2;c)$ to $(x_0,x_1,x_2;-c)$.
So the birational map is the simplest example of flop
between 3 dimensional varieties and hence
this is why $I$ is not an isomorphism
(cf. \cite{Utr} and \cite{Kol}, 12.1. Example). 
Similarly the cycle $C_4-C_3$ is the vanishing cycle around $c=-1$.
\section{Regular functions on $\x[t_0,c_0] \backslash D$.}
\setcounter{equation}{0}
\par 
We begin with a lemma.
\begin{Lemma}
If there exists a complete curve $C$ on $\x[t_0,c_0] \backslash \D[t_0,c_0]$,
then $C \cap W_1 \subset W_1$ is an open algebraic curve.
\end{Lemma}

\Proof
It follows from the construction of the surface $\x[t_0,c_0]$
that $\x[t_0,c_0] \backslash W_1= \D[t_0,c_0] \cup C_1 \cup C_4$.
The curve $C_1$ intersects $\D[t_0,c_0]$ on $D_1$ and the curve $C_4$ intersects
$\D[t_0,c_0]$ on $D_7$.
\begin{Proposition}
{\bf (i)} If the parameter $c_0$ is not an integer, then there is no 
projective algebraic curve on $\x[t_0,c_0] \backslash \D[t_0,c_0]$.
{\bf (ii)} If $c_0$ is an integer, then there is only one projective algebraic
curve $C$ on $\x[t_0,c_0] \backslash \D[t_0,c_0]$.
The curve $C$ is isomorphic to $\P^1$ and $C^2=-2$.
In other words, $C$ is a $-2$-curve.
\end{Proposition}

\Proof
Let us assume that there exists a complete curve $C$ on 
$\x[t_0,c_0] \backslash \D[t_0,c_0]$.
Since the Painlev\'e equations have no movable singular points and since 
$\x[t_0,c_0] \backslash \D[t_0,c_0]$ is the space of initial conditions, 
integrating the curve $C$ along the foliation on 
$\x[t_0,c_0] \backslash \D[t_0,c_0]$ generated by the Hamiltonian system
\begin{eqnarray*}
S_2(c_0) 
\left\{
\begin{array}{l}
\frac{\displaystyle dy_1}{\displaystyle dt} =
y_1^2+z_1+\frac{\displaystyle t}{\displaystyle 2}, \\
\smallskip \\
\frac{\displaystyle dz_1}{\displaystyle dt}=- 2 y_1 z_1 + c_0
\end{array}
\right.
\end{eqnarray*}
on $W_1$, we get an analytic family  
${\frak C} \rightarrow \left( \Spec \C[t] \right)^{an}$
of complete algebraic curves on $\x[t,c_0] \rightarrow \Spec \C[t]$.
So by G.A.G.A., there exists a field extension $K \supset \C(t)$ 
such that $K$ is a subfield of the field of meromorphic functions
on an open set of $\C$ with coordinate system $t$ and 
such that the analytic family
${\frak C} \rightarrow \left( \Spec \C[t] \right)^{an}$
is an algebraic family of curves defined over the field $K$.
We may assume that the field is closed under the derivation $d/dt$.
Let $F(y_1,z_1)=0$ be a defining equation of 
$\left( W_1 \times \C \cap {\frak C} \right) \rightarrow \Spec K$ on 
$\left( \x[t,c_0] \cap W_1 \right) \otimes_{\C[t]} K \rightarrow \Spec K$
 so that $F(y_1,z_1) \in K[y_1,z_1]$.
Since $W_1 \times \C \cap {\frak C} $ is invariant under the Hamiltonian flow
\begin{eqnarray*}
D(c_0)=\frac{\partial}{\partial t} + \left( y_1^2+z_1 +\frac{t}{2} \right)
\frac{\partial}{\partial y_1} +\left( c_0 -2y_1 z_1 \right) 
\frac{\partial}{\partial z_1},
\end{eqnarray*}
$f(y_1,z_1)$ divides $D(c)\left( f(y_1,z_1) \right)$ in the polynomial ring 
$K[y_1,z_1]$.
In the language of \cite{Uir}, $F$ is an invariant polynomial over $K$.
Then Theorem (2.1), \cite{UW} 
implies $c_0 = \alpha_0 -\frac{1}{2} \in {\bf Z}$
and the curve $C$ arises from the Riccati equation and $C^2=-2$.
\begin{Lemma}
The canonical divisor $K$ of $\x[t_0,c_0]$ is given by a divisor $-F$,
where
\setcounter{meq}{\value{equation}}
\setcounter{equation}{0}
\renewcommand{\theequation}{\thenewsection.\themeq.\arabic{equation}}
\begin{eqnarray}
F=2D_0+D_1+2D_2+3D_3+4D_4+3D_5+2D_6+D_7.
\end{eqnarray}
\setcounter{equation}{\value{meq}}
\end{Lemma}

\Proof
The lemma follows from the following two observations
and the construction of $\x[t_0,c_0]$. 
First, the canonical divisor $K_0$ of $Z_0[t_0,c_0]$ is given by
$-2S$.
Second, let $K_i$ be canonical divisor of $Z_i[t_0,c_0]$ $( 0 \leq i \leq 8)$.
Then $K_{i+1} \sim \pi_{i+1}^{*} K_i+E_{i+1}$, where $E_{i}$ is exceptional
divisor of $\pi_{i+1}:Z_{i+1}[t_0,c_0] \rightarrow Z_i[t_0,c_0]$ 
$0 \leq i \leq 7$.\\
\par
$F$ is the null vector of the Cartan matrix $C$ of $\widetilde{A}_1$.
Namely $C {}^{t}(2,1,2,3,4,3,2,1)=0$.
\begin{Proposition}
$h^0 \left( \x[t_0,c_0],-K \right)=1$.
\end{Proposition}

\Proof 
Since $F$ is effective, $h^0 \left( \x[t_0,c_0],-K\right) \geq 1$.
Assume $h^0 \left( \x[t_0,c_0],-K\right) \geq 2$.
Then there exists an effective divisor $H$ linearly equivalent to $-K$
and distinct of $F$.
Let $H=H_1+H_2$ be a decomposition into two effective divisors
such that the support of $H_1$ is a subset of
 $\{D_i \, | \,   0 \leq i\leq 7\}$
and such that no irreducible component of $H_2$ is $D_i$, $0 \leq i \leq 7$.
We show $H_2=0$. Let us assume $H_2 \neq 0$.
Since we blow up the ruled surface $Z_0[t_0,c_0]$ 8 times,
\begin{eqnarray}
0=K^2=(K\,.\,-H)=(K\,.\,-H_1-H_2)=-(K\,.\,H_1)-(K\,.\,H_2).
\end{eqnarray}
Since $D_i$ is a $-2$-curve, $(K\,.\,D_i)=0$ for $0 \leq i \leq 7$.
So 
\begin{eqnarray}
(K\,.\,H_1)=0.
\end{eqnarray}
First assume that $c_0$ is not an integer.
Then there is no projective curve on 
$\x[t_0,c_0] \backslash \bigcup_{i=0}^{7} D_i$
by Lemma (4.1).
Hence $(-K\,.\,H_2)=(F\,.\,H_2) >0$.
This contradicts (4.5) and (4.6).
Therefore $H_2=0$ or the support of $H \subset \bigcup^{7}_{i=0} D_i$.
Since the $D_i$'s are linearly
independent divisors in $\Pic \x[t_0,c_0]$,
this shows $H=F$, which is contradiction.
If $c_0$ is an integer, there exists the unique projective curve
$C$ on $\x[t_0,c_0] \backslash \bigcup_{i=0}^7 D_i$.
If we assume $H_2 \neq 0$, then 
the above argument shows $H_2=nC$ for an appropriate integer $n \geq 0$.
Since $C$ is a $-2$-curve,
\begin{eqnarray*}
0=(K.C)=(-H.C)=(-H_1-H_2.C)=(-H_1.C)-(H_2.C)=-n(C.C)=2n.
\end{eqnarray*}
So we conclude $H_1=0$ and the argument above leads us to a contradiction. 
\begin{Corollary}[to the proof.]
For every integer $m \geq 0$,
\begin{eqnarray*}
H^0(\x[t_0,c_0], -mK)=\C.
\end{eqnarray*}
\end{Corollary}

\Proof The assertion being trivial for $m=0$, we may assume $m>0$. 
Then the argument of the proof of Proposition (4.4)
allows us to prove the corollary.
We can formulate Corollary (4.7) in another form.
\begin{Corollary}
\begin{eqnarray*}
H^0(\x[t_0,c_0] \backslash \D[t_0,c_0], {\cal O})=\C.
\end{eqnarray*}
\end{Corollary}

\Proof
Since $-K \sim 2D_0+D_1+2D_2+3D_3+4D_4+3D_5+2D_6+D_7$,
this result follows from Corollary (4.7).
\section{Takano coordinate systems.}
\setcounter{equation}{0}
\par
It follows from Proposition (4.4) that there exists a $2$-form $\omega$
on $\x[t_0,c_0] \backslash \D[t_0,c_0]$ such that $\omega$ vanishes
at no point of $\x[t_0,c_0] \backslash \D[t_0,c_0]$.
Moreover $\omega$ is  unique up to a non-zero constant multiplication.
Namely, there exists a symplectic structure on 
$\x[t_0,c_0] \backslash \D[t_0,c_0]$.
K. Takano and his colleagues introduced nice coordinate systems on
$\x[t_0,c_0] \backslash \D[t_0,c_0]$
such that locally the $2$-form is written in a canonical form in terms of the 
coordinate system.
They cover $\x[t_0,c_0] \backslash \D[t_0,c_0]$ by there copies of $\C^2$
such that coordinate transformations are symplectic.
We call these coordinate systems the Takano coordinate systems.
To explain the Takano coordinate systems, let us come back to the construction
of $\x[t_0,c_0] \backslash \D[t_0,c_0]$.
To construct $\x[t_0,c_0]$, we started from the ruled surface $Z_0[t_0,c_0]$.
We blow up $Z_0[t_0,c_0]$ eight times and removed the proper transform 
$D_0$ of $S$,
the proper transforms $D_1, D_2, \ldots, D_7$ of all the exceptional divisor 
except for the last exceptional divisor $E_8=C_3=D_8$.
Let us recall
\begin{eqnarray*}
S=
\left\{
(y_2,z_2) \in W_2 \, | \,   z_2=0
\right\}
\cup
\left\{
(y_4,z_4) \in W_4 \, | \,   z_4=0
\right\}
\end{eqnarray*}
so that 
\begin{eqnarray*}
Z_0[t_0,c_0]=W_1 \cup W_3 \cup S.
\end{eqnarray*}
Hence
\begin{eqnarray}
\x[t_0,c_0] \backslash \D[t_0,c_0]
=
W_1 \cup W_3 \cup 
\left( 
D_8 \backslash (D_8 \cap D_7)
\right)
\end{eqnarray}
(cf. Fig.(1.3)).
So we start from $W_1 \cup W_3 \cup W_4$.
We use the notation of \S2, \cite{Utr}.
We have to resolve the rational map 
\begin{eqnarray}
F:W_4 \cdots \rightarrow \P^1,
\hspace{0.5cm}
(y_4,z_4) \mapsto (y_4^4z_4, 2z_4-y_4^4+ty_4^2 z_4+(2c+1)y_4^3z_4)
\end{eqnarray}
(cf. Sublemma (3.6) in \cite{Utr}).
Let $(y_{4+i}, z_{4+i})$ be the coordinate system of $W_{4+i}(z)=\C^2$
$1 \leq i \leq 4$ so that we have $y_{4+i}=y_4$, 
$z_{4+i-1} = y_{4+i} z_{4+i}$ $(1 \leq i \leq 4)$.
Hence
\begin{eqnarray}
y_4=y_5=y_6=y_7=y_8, \hspace{0.5cm} y_4=y_8^4z_8.
\end{eqnarray} 
On $W_8(z)$, $F$ is written as
\begin{eqnarray}
F:W_8(z) \rightarrow \P^1,
\hspace{0.5cm}
(y_8^4z_8, 2 z_8 -1+t y_8^2 z_8 + (2c+1) y_8^3).
\end{eqnarray}
Now we introduce a coordinate system 
$(y_8,v_8)$ 
on $W_8(z)$
such that 
\begin{eqnarray}
v_8=\frac{1}{z_8}.
\end{eqnarray}
Precisely speaking, we consider $W_8'(z)=\A^2$ with
coordinate system $(y_8,z_8)$ and birational map
\begin{eqnarray}
W_8'(z) \cdots \rightarrow W_8(z), \hspace{0.5cm}
(y_8,v_8) \mapsto (y_8,1/v_8).
\end{eqnarray}
Let $W_8'(z)^{0}=\{(y_8,v_8) \in W_9'(z)\, | \,  v_8 \neq 0\}$
so that $W_8'(z)^{0}$ is identified with an open set of $W_8(z)$ by (5.6).
We have on $W_8'(z)$
\begin{eqnarray}
F:W_8'(z) \rightarrow \P^1, \hspace{0.5cm} 
(y_8,v_8) \mapsto (y_8^4,2-v_8+ty_8^2+(2c+1)y_8^3).
\end{eqnarray}
So $(y_8,v_8)=(0,2)$ is the singular point of $F$ in (5.7).
We blow up $W_8'(z)$ at $(y_8,v_8)=(0,2)$ 
to resolve the rational map $F$ given by (5.7).
We set
\begin{eqnarray}
v_8-2=y_9z_9, \hspace{0.5cm} y_8=y_9.
\end{eqnarray} 
On $W_9=\A^2$ with coordinate system $(y_9,z_9)$, we have 
\begin{eqnarray}
F:W_9 \rightarrow \P^1, \hspace{0.5cm}
(y_9,z_9) \mapsto (y_9^3,-z_9+ty_9+(2c+1)y_9^2).
\end{eqnarray}
The singular point of $F$ in (5.9) is $(y_9,z_9)=(0,0)$.
We blow up $W_9$ at this point.
So we set 
\begin{eqnarray}
y_{10}=y_9, \hspace{0.5cm} z_9=y_9 z_{10}.
\end{eqnarray}
On $W_{10}=\A^2$ with coordinate system $(y_{10},z_{10})$, we have 
\begin{eqnarray}
F:W_{10} \rightarrow \P^1, \hspace{0.5cm}
(y_{10}, z_{10}) \mapsto (y_{10}^2,-z_{10}+t+(2c+1)y_{10}).
\end{eqnarray}
The singular point of $F$ in (5.11) is $(y_{10},z_{10})=(0,t)$.
We blow up $W_{10}$ at (0,t).
So we set 
\begin{eqnarray}
y_{11}=y_{10}, \hspace{0.5cm}z_{10}-t=y_{11}z_{11}.
\end{eqnarray}
On $W_{11}=\A^2$ the coordinate system $(y_{11},z_{11})$,
we have 
\begin{eqnarray}
W_{11} \rightarrow \P^1, \hspace{0.5cm} 
(y_{11},z_{11}) \mapsto (y_{11}, -z_{11}+(2c+1)).
\end{eqnarray}
The singular point of (5.11) is $(y_{11},z_{11})=(0,2c+1).$
We blow up $W_{11}$ at $(0,2c+1)$.
So we set 
\begin{eqnarray}
y_{12}=y_{11}, \hspace{0.5cm} z_{11}-(2c+1) = y_{12}z_{12}.
\end{eqnarray}
On $W_{12}=\C^2$ with coordinate system $(y_{12},z_{12})$,
we have 
\begin{eqnarray}
F:W_{12} \rightarrow \P^1, \hspace{0.5cm} (y_{12},z_{12}) \mapsto (1,z_{12}).
\end{eqnarray}
Namely the rational map is resolved on $W_{12}$.
It follows from (5.8), (5.10), (5.12), (5.14)
\begin{eqnarray}
v_8=y_{12}^4z_{12}+(2c+1)y_{12}^3+ty_{12}^2+2.
\end{eqnarray}
So we have a rational map 
\begin{eqnarray}
W_{12} \cdots \rightarrow W_8'(z) \rightarrow W_8(z)
\end{eqnarray}
sending $(y_{12},z_{12})$ to 
\begin{eqnarray*}
(y_8,z_8)=
\left( y_{12}, \frac{1}{y_{12}^4z_{12}+(2c+1)y_{12}^3+ty_{12}^2+2}\right).
\end{eqnarray*}
We have to restrict the map (5.17) on an open set
\begin{eqnarray*}
W_{12}^{0}=
\left\{
(y_{12},z_{12})\, | \,  
y_{12}^4z_{12}+(2c+1)y_{12}^3+ty_{12}^2+2\neq 0
\right\}.
\end{eqnarray*}
Then we have a regular map 
\begin{eqnarray}
&W_{12}^{0} \rightarrow W_8(z),&\\
\smallskip \nonumber \\
&(y_{12},z_{12}) \mapsto (y_8,z_8)=
\left( y_{12}, \frac
{\displaystyle 1}
{\displaystyle y_{12}^4z_{12}+(2c+1)y_{12}^3+ty_{12}^2+2}\right).
&
\nonumber
\end{eqnarray}
Hence it follows from (5.3), (5,18) that we have a regular map
\begin{eqnarray}
&W_{12}^{0} \rightarrow W_4(z),&\\
\smallskip \nonumber \\
&(y_{12},z_{12}) \mapsto (y_8,z_8)=
\left( y_{12}, \frac
{\displaystyle y_{12}^4}
{\displaystyle y_{12}^4z_{12}+(2c+1)y_{12}^3+ty_{12}^2+2}\right).
&
\nonumber
\end{eqnarray}
Since $z_4=\frac{1}{z_3}$, we have by (5.19)
\begin{eqnarray*}
\frac{1}{z_3}=
\frac{y_{12}^4}{y_{12}^4z_{12}+(2c+1)y_{12}^3+ty_{12}^2+2}
\end{eqnarray*}
or 
\begin{eqnarray}
z_3=z_{12}+\frac{2c+1}{y_{12}}+\frac{t}{y_{12}^2}+\frac{2}{y_{12}^4}.
\end{eqnarray}
So $\x[t_0,c_0] \backslash \D[t_0,c_0]$ is covered by $W_1$, $W_3$ and 
$W_{12}^{0}$.
A point $(y_3,z_3) \in W_3$ and a point $(y_{12},z_{12}) \in W_{12}^{0}$
are identified if {\bf (i)} $y_4=y_{12} \neq 0$ 
and if {\bf (ii)} we have (5.20).
In view of (5.18), we may enlarge $W_{12}^{0}$ or 
we may replace it by $W_{12}$.
So in view of (5.1), we have proved the following

\begin{Theorem}
$\x[t_0,c_0] \backslash \D[t_0,c_0]$ is covered by
three copies $W_1$, $W_3$, $W_{12}$ of $\A^2$.
We glue together these copies by following rule.\\
{\bf (i)} A point $(y_1,z_1) \in W_1$ and a point $(y_3,z_3)$ are identified if
$y_1y_3=1$ and if $z_1=y_3(c-y_3 z_3)$.\\
{\bf (ii)}
A point $(y_1,z_1) \in W_1$ and a point $(y_{12},z_{12}) \in W_{12} $
are identified if $y_1y_{12}=1$ and 
if $z_1+2y_1^2+t_0 =y_{12}\left(-(c+1)-y_{12} z_{12} \right)$.\\
{\bf (iii)}
 A point $(y_3,z_3) \in W_3$ and a point $(y_{12},z_{12}) \in W_{12}$
are identified if $y_3=y_{12} \neq 0$ and if 
we have 
\begin{eqnarray*}
z_3=z_{12}+\frac{2c+1}{y_{12}}+\frac{t}{y_{12}^2}+\frac{2}{y_{12}^4}.
\end{eqnarray*}
\end{Theorem}

\Proof
It is sufficient to see that {(ii)} is a consequence of {(i)}
and {(iii)}.
We can check this by an easy calculation.
\begin{Corollary}
$2$-forms $dy_1 \wedge dz_1$ on $W_1$, $dy_3 \wedge dz_3$ on $W_3$
and $dy_{12} \wedge dz_{12}$ on $W_{12}$
coincide on the intersections $W_i \cap W_j$
and thus define a symplectic structure on $\x[t_0,c_0] \backslash \D[t_0,c_0]$.
\end{Corollary}

\Proof
This is an immediate consequence of Theorem (5.21).
\begin{Remark}
If we notice the relation {(iii)} follows from {(i)} and {(ii)},
we have a symmetry
\begin{eqnarray*}
y_1 \mapsto -y_1,
\hspace{0.5cm}
z_1 \mapsto -(z_1+2y_1^2+t),
\hspace{0.5cm}
c \mapsto -(c+1),
\end{eqnarray*}
which sends as a consequence $z_1+2y_1^2+t \mapsto -z_1$,
$-(c+1) \mapsto -c$, of 
$W_1 \cup W_3 \cup W_{12}=\x[t_0,c_0] \backslash \D[t_0,c_0]$. 
\end{Remark}

The open surface $\x[t_0,c_0] \backslash D$ is covered by
three copies $W_1$, $W_3$, $W_{12}$ of $\A^2$.
As Corollary (5.23) shows, the $2$-forms
$dy_1 \wedge dz_1$, 
$dy_3 \wedge dz_3$, 
$dy_{12} \wedge dz_{12}$
glue together and define a $2$-form $\omega$ on $\x[t_0,c_0] \backslash D$.
The vanishing cycles span 
$L\left( A_1^{(1)}\right) = L\left(E_7^{(1)} \right)^{\perp}$
and hence
by (2.30)
${\bf Z} (C_2-C_1) \oplus {\bf Z}(C_4-C_3)
=H_2(\x[t_0,c_0] \backslash D, {\bf Z})
$.
\begin{Proposition}
\begin{eqnarray*}
\int_{C_2-C_1} \omega =c, \hspace{0.5cm} \int_{C_4-C_3}\omega =-c-1.
\end{eqnarray*}
\end{Proposition}

\Proof
We show the first equality.
If $c=0$, $C_2-C_1$ is a complete curve on $\x[t_0,c_0] \backslash D$ so that 
\begin{eqnarray*}
\int_{C_2-C_1} \omega =0
\end{eqnarray*}
because $\omega$ is a holomorphic 2-form.
So we may assume $c \neq 0$.
Since $y_3=y_4$, $z_3=1/z_4$, we have on $W_4$ 
\begin{eqnarray*}
\omega=dy_3 \wedge dz_3 =-dy_4 \wedge \frac{1}{z_4^2}dz_4.
\end{eqnarray*}
Blowing up $W_4$ at $(y_4,z_4)=(0,0)$,
we get a morphism $\widetilde{W}_4 \rightarrow W_4$.
The surface $\widetilde{W}_4$ is covered by $W_4(y)$ and $W_4(z)$
that are isomorphic to $\A^2$.
We have coordinate systems $(Y,z)$ on $W_4(y) \simeq \A^2$ and
$(y,Z)$ on $W_4(z) \simeq \A^2$.
We identify a point $(Y,z) \in W_4(y)$ and $(y,Z) \in W_4(z)$
if $y=Yz$, $z=yZ$ and if $YZ=1$.
The morphism $\widetilde{W}_4=W_4(y) \cap W_4(z) \to W_4$ is defined respectively on $W_4(y)$ by
$$
W_4(y) \rightarrow W_4=\A^2, \qquad 
(Y,z) \mapsto (Yz,z)
$$
and on $W_4(z)$ by
$$
W_4(y) \rightarrow W_4=\A^2, \qquad  
(y,Z) \mapsto (y,yZ).
$$
The curves $C_1$ and $C_2$ are defined on $W_4$ respectively by
\begin{eqnarray*}
C_1 \cap W_4 = \left\{ (y_4,z_4) \in \A^2 \, | \,   y_4=0\right\}
\end{eqnarray*}
and 
\begin{eqnarray*}
C_2 \cap W_4 = \left\{ (y_4,z_4) \in \A^2 \, | \,   y_4-cz_4=0\right\}.
\end{eqnarray*}
So we have on $W_4(y) \simeq \A^2$,
\begin{eqnarray*}
C_1 \cap W_4(y) = \left\{ (Y,z) \in W_4(y) \, | \,   Y=0\right\}
\end{eqnarray*}
and 
\begin{eqnarray*}
C_2 \cap W_4(y) = \left\{ (Y,z) \in W_4(y) \, | \,   Y=c\right\}.
\end{eqnarray*}
Now $D_1 \cap W_1(y) =\left\{(Y,z) \in  \A^2 \, | \,  z=0
\right\}$
and hence on $W_4(y)$, $C_1 \cap D_1= \left\{ (0,0)\right\}$,
$C_2 \cap D_1 = \left\{ (c,0)\right\}$.
Now let $\gamma$ be a segment in 
$D_1 \cap W_4(y)= \left\{ (Y,z) \in W_4(y) \, | \,   z=0 \right\}$
joining the points $(c,0)$ and $(0,0)$.
Let $\tau$ be a closed tublar neighborhood of $\gamma$.
We set $C_1 \cap \tau =\tau_1$, $C_2 \cap \tau =\tau_2$.
We have on $W_4(y)$ 
\begin{eqnarray*}
-dy_4 \wedge \frac{1}{z_4^2}dz_4=-dY \wedge \frac{dz}{z}.
\end{eqnarray*} 
Since $\partial \tau$ is homologous to $C_2-C_1$ in $\x[t_0,c_0] \backslash D$,
\begin{eqnarray*}
\int_{C_2-C_1} \omega =\int_{\partial \tau} \omega
=\int_{\partial \tau} -dY \wedge \frac{dz}{z}
=\int_{\gamma} -2 \pi i dY
=2\pi i c.
\end{eqnarray*}
The isomorphism
\begin{eqnarray*}
J(c_0,-c_0-1,t_0):\x[t_0,c_0] \rightarrow \x[t_0,-c_0-1]
\end{eqnarray*}
maps $C_3[t_0,c_0]$ to $C_1[t_0,-c_0-1]$, $C_4[t_0,c_0]$ to
$C_2[t_0,-c_0-1]$ and $\omega[t_0,c_0]$ to $\omega[t_0,-c_0-1]$.
So
\begin{eqnarray*}
\int_{C_4[t_0,c_0]-C_3[t_0,c_0]} \omega[t_0,c_0]
=
\int_{C_2[t_0,-c_0-1]-C_1[t_0,-c_0-1]} \omega[t_0,-c_0-1]
=-c_0-1.
\end{eqnarray*}
\vskip 5mm

In particular, since either $c \neq 0$ or $-c-1 \neq 0$,
we have proved the following 
\begin{Corollary}
The de Rham class of $\omega$ in 
$H^2\left( \x[t_0,c_0] \backslash D, \C \right)$ is not 0.
In other words the closed 2-form $\omega$ on $x[t_0,c_0] \backslash D$
is not exact.
\end{Corollary}
\begin{Theorem} {\rm (Cf. \cite{A}). }
We have 
\begin{eqnarray*}
H^0(\x [t_0,c_0] \backslash D, \Theta)=0
\end{eqnarray*}
or equivalently
\begin{eqnarray*}
H^0(\x[t_0,c_0] \backslash D, \Omega^1)=0.
\end{eqnarray*}
\end{Theorem}

\Proof 
Since we have a non-degenerate 2-form $\omega$ on $\x [t_0,c_0] \backslash D$,
the sheaf $\Theta$ is isomorphic to $\Omega^1$.
So we have to prove $H^0 \left(\x [t_0,c_0] \backslash D, \Omega^1\right)=0$.
Let $\eta \in H^0 \left(\x[t_0,c_0] \backslash D, \Omega^1\right)$.
We show $\eta=0$.
In fact, $d\eta \in H^0 \left(\x[t_0,c_0] \backslash D, \Omega^2\right)$.
It follows from Proposition (4.4)
that we can find a complex number $\lambda$ such that $d\eta=\lambda \omega$.
By Corollary (5.25), $\lambda=0$ or $d \eta=0$.
The algebraic surface $\x [t_0,c_0] \backslash D$ is covered by
three copies $W_1$, $W_2$, $W_{12}$ of $\A^2$ so that we can find
$\varphi_i \in H^0(W_i, {\cal O})$ for $i=1,3,12$
such that $d \varphi_i = \eta$ on $W_i$ for $i=1,3,12$.
Hence $\{\varphi_i-\varphi_j\}_{W_i\cap W_j}$ is a 1-cocycle
or is an element of $Z^1(\bigcup W_i, \C)$.
Since $H_{zar}^1\left( \x [t_0,c_0] \backslash D, \C \right) =0$,
we can find constants $k_i \in \C$ for $i=1,3,12$ such that
\begin{eqnarray*}
\varphi_i-\varphi_j=k_i-k_j.
\end{eqnarray*} 
In other words the functions $\varphi_i-k_i$ glue together to
give an element of $H^0\left( \x [t_0,c_0] \backslash D, {\cal O} \right)$.
Now it follows from Corollary (4.8) that this function is constant $k$.
Therefore we have $\varphi_i-k_i=k$ on $W_i$.
This shows $\varphi_i=k+k_i \in \C $.
Hence $0=d\varphi_i=\eta$.
This is what we had to show.
\section{Hamiltonian system.}
\setcounter{equation}{0}

\par
We worked in \S5 on $\x[t_0,c_0]$ for a fixed $t_0$ and $c_0$.
We can apply this argument to the relative case $\x \slash \Spec \C[t,c]$.
In theory of Painlev\'e equations, however we have to study 
$\x \slash \Spec \C[c]$.
We have three copies of $\C^4$:
\begin{eqnarray*}
&&{\cal W}_1:=\Spec \C[y_1,z_1,t,c],\\
\smallskip \nonumber\\
&&{\cal W}_3:=\Spec \C[y_3,z_3,t,c],\\
\smallskip \nonumber\\
&&{\cal W}_{12}:=\Spec \C[y_{12},z_{12},t,c]. 
\end{eqnarray*}  
A point $(y_1,z_1,t_1,c_1) \in {\cal W}_1$ and a point 
$(y_3,z_3,t_3,c_3) \in {\cal W}_3$
and identified if $(t_1,c_1)=(t_3,c_3)$, $y_1 y_3=1$ 
and if $z_1=y_3(c-y_3 z_3)$.
A point $(y_1,z_1,t_1,c_1)$ of ${\cal W}_1$ and a point 
$(y_{12},z_{12}, t_{12}, c_{12})$ of ${\cal W}_{12}$
are identified if $(t_1,c_1)=(t_{12},c_{12})$, $y_1 y_{12}=1$ and if
\begin{eqnarray*}
z_1+2y_1^2+t_1=y_{12}\left( -(c+1) -y_{12}z_{12}\right).
\end{eqnarray*}
A point $(y_3,z_3,t_3,c_3)$ of ${\cal W}_3$ and a point 
$(y_{12},z_{12}, t_{12}, c_{12})$ of ${\cal W}_{12}$ 
are identified if $(t_3,c_3)=(t_{12},c_{12})$ and if
\begin{eqnarray*}
z_3=z_{12}+\frac{2c_{12}+1}{y_{12}}+\frac{t_{12}}{y_{12}^2}+\frac{2}{y_{12}^4}.
\end{eqnarray*}
By gluing ${\cal W}_1$, ${\cal W}_2$ and ${\cal W}_3$ together by this rule, 
we have 
${\cal W}_1 \cup {\cal W}_2 \cup {\cal W}_3=\x \backslash \D$.
Now we consider $\x$ as a variety over $\Spec \C[c]$.
In other words $\x \rightarrow \Spec \C[c]$ is a family of threefolds
parameterized by $\Spec \C[c]$.
The differential, as well the $\wedge$, is take over $\C[c]$.
To distinguish the differential over $\C[c]$ from the one over $\C[t,c]$,
we denote the former by $d_{/c}$ so that $d_{/c} \C[c]=0$ but 
$d_{/c} t \neq 0$. 
Similarly we denote by $\wedge_{/c}$ the wedge over $\C[c]$.

We look for polynomials $H_1(t,c,y_1,z_1) \in \C[t,c,y_1,z_1]$,
$H_3(t,c,y_3,z_3) \in \C[t,c,y_3,z_3]$, 
$H_{12}(t,c,y_{12}, z_{12}) \in \C[t,c,y_{12},z_{12}]$
such that
\begin{eqnarray}
d_{/c}y_1 \wedge_{/c} d_{/c}z_1 + d_{/c} H_1 \wedge_{/c} d_{/c}t,\\
\smallskip \nonumber \\
d_{/c}y_3 \wedge_{/c} d_{/c}z_3 + d_{/c} H_3 \wedge_{/c} d_{/c}t \nonumber
\end{eqnarray}
and 
\begin{eqnarray}
d_{/c}y_{12} \wedge_{/c} d_{/c}z_{12} + d_{/c} H_{12} \wedge_{/c} d_{/c}t
\nonumber
\end{eqnarray}
glue together and 
define a $2$-form on $\x \backslash \D$ over $\Spec \C[c]$.
\begin{Lemma}
Let $H_1=y_1^2 z_1 +\frac{1}{2} z_1^2+\frac{t}{2} z_1 - c y_1 
\in \C[t,c,y_1,z_1]$.
Then, \\
{\bf (i)} $H_3(t,c,y_3,z_3):=H_1(t,c,\frac{1}{y_3}, y_3(c-y_3z_3))$
is a polynomial in $t$, $c$, $y_3$, $z_3$.\\
{\bf (ii)} $H_{12}(t,c,y_{12},z_{12}):=
H_3 \left( t,c,y_{12}, 
z_{12}+\frac{2c+1}{y_{12}}+\frac{t}{y_{12}^2}+\frac{2}{y_{12}^4} \right)
-\frac{1}{y_{12}}$ is a polynomial in $t$, $c$, $y_{12}$, $z_{12}$.
\end{Lemma}

\Proof
The lemma is proved by a simple calculation.
\begin{Lemma}
We have
\setcounter{meq}{\value{equation}}
\setcounter{equation}{0}
{
\renewcommand{\theequation}{\thenewsection.\themeq.\arabic{equation}}
\begin{eqnarray}
d_{/c} y_1 \wedge_{/c} d_{/c} z_1 = d_{/c} y_3 \wedge_{/c} d_{/c} z_3
\end{eqnarray}
and 
\begin{eqnarray}
d_{/c} y_3 \wedge_{/c} d_{/c} z_3 =d_{/c} y_{12} \wedge_{/c} d_{/c} z_{12}
-\frac{1}{y_{12}^2} d_{/c} y_{12} \wedge_{/c} d_{/c} t.
\end{eqnarray}
}
\setcounter{equation}{\value{meq}}
\end{Lemma}

\Proof These formulas are consequences of the identification rule of 
${\cal W}_1$, ${\cal W}_2$ and ${\cal W}_{12}$.
\begin{Proposition}
Closed 2-forms 
$d_{/c} y_1 \wedge_{/c} d_{/c} z_1 +
d_{/c} H_1 \wedge_{/c} d_{/c} t$, 
$d_{/c} y_3 \wedge_{/c} d_{/c} z_3
+d_{/c} H_3 \wedge_{/c} d_{/c} t$,
$d_{/c} y_{12} \wedge_{/c} d_{/c} z_{12} +
d_{/c} H_{12} \wedge_{/c} d_{/c} t$
glue together and define a 2-form $\omega_c$ on 
$\x \backslash \D$ over $\Spec \C[c]$.
\end{Proposition}

\Proof 
We have to show 
\addtocounter{equation}{1}
\setcounter{meq}{\value{equation}}
\setcounter{equation}{0}
{
\renewcommand{\theequation}{\thenewsection.\themeq.\arabic{equation}}
\begin{eqnarray}
d_{/c} y_1 \wedge_{/c} d_{/c} z_1 +
d_{/c} H_1 \wedge_{/c} d_{/c} t
=
d_{/c} y_3 \wedge_{/c} d_{/c} z_3 +
d_{/c} H_3 \wedge_{/c} d_{/c} t
\end{eqnarray}
on ${\cal W}_1 \cap {\cal W}_3$
and 
\begin{eqnarray}
d_{/c} y_3 \wedge_{/c} d_{/c} z_3 +
d_{/c} H_3 \wedge_{/c} d_{/c} t
=
d_{/c} y_{12} \wedge_{/c} d_{/c} z_{12} +
d_{/c} H_{12} \wedge_{/c} d_{/c} t
\end{eqnarray}
}
\setcounter{equation}{\value{meq}}
on ${\cal W}_3 \cap {\cal W}_{12}$.
(6.5.1) follows from (6.3.1) and the fact that $H_1=H_3$ 
(cf. Lemma (6.2), {(ii)}). 
On the other hand (6.5.2) follows from Lemma (6.2), {(ii)} and (6.3.2).
\par 
\vspace{1pc}
The 2-form $\omega_c$ defines a skew symmetric form on 
$\Theta_{\x \backslash \D \slash \Spec \C[c]}$, of which the null foliation
is the second Painlev\'e equation.
Namely $\Theta_{\x \backslash \D \slash \Spec \C[c]}$ is a vector bundle
of rank three on which we have the skew symmetric form $\omega_c$.
So at each point $p$ of $\x \backslash \D$, 
there exist a 1-dimensional subspace $\Sigma_p $ of 
$\Theta_{\x \backslash \D \slash \Spec \C[t]}$ such that
$\omega_c(\Sigma_p, \Theta_{\x \backslash \D \slash \Spec \C}, p)=0$.
The subspace $\Sigma_p$ defines a foliation.
Rigorously speaking, we have to fix $c=c_0 \in \C$ and work on the fiber
$\x_{c_0} \backslash \D_{c_0}$, which is a threefold, over $c=c_0$ of 
$\x \backslash \D \rightarrow \Spec \C[c]$.
Locally on ${\cal W}_1$ over $\C[c]$ for example the null foliation is given by
the Hamiltonian flow
\begin{eqnarray}
\frac{y_1}{dt}=\frac{\partial H_1}{ \partial z_1}, \hspace{0.5cm}
\frac{z_1}{dt}=-\frac{\partial H_1}{\partial y_1}.
\end{eqnarray} 
Namely
\begin{eqnarray*}
\left\{
\begin{array}{l}
\frac{\displaystyle dy_1}{\displaystyle dt}=y_1^2+z_1
+\frac{\displaystyle t}{\displaystyle 2},\\
\smallskip \nonumber\\
\frac{\displaystyle dz_1}{\displaystyle dt}=-2y_1z_1+c.
\end{array}
\right.
\end{eqnarray*}
\begin{Theorem}
The Hamiltonian functions $H_1$, $H_3$, $H_{12}$ are unique. Namely
if $H'_1$, $H'_3$, $H'_{12}$ are polynomials such that the 
2-forms in (6.1) glue together, then $H_i-H'_i$ does not depend on $y_i$, $z_i$
for $i=1,3,12$. 
\end{Theorem}

\Proof
Let $H'_i \in \C[t,c,y_1,z_1]$, $H'_3 \in \C[t,c,y_3,z_3]$, 
$H'_{12} \in \C[t,c,y_{12}, z_{12}]$ such that the 2-forms 
of (6.1) glue together. 
We have to show that $H_1-H'_1$, $H_3-H'_3$, 
$H_{12}-H'_{12}$ are functions of $t$, $c$.
It follows from (6.1), we have $d(H_1-H'_1)=d(H_3-H'_3)$ on 
${\cal W}_1 \cap {\cal W}_3$,
$d(H_1-H'_1)=d(H_{12}-H'_{12})$ on 
${\cal W}_1 \cap {\cal W}_{12}$,
where $d$ is taken over $\C[t,c]$ so that $d(\C[t,c])=0$.
We set
\begin{eqnarray*}
\xi_{ij}=(H_i-H'_i)-(H_j-H'_j)
\end{eqnarray*}
for $i,j \in \{1,3,12\}$.
Then $\xi_{i,j} \in \C$ and $\{\xi_{ij}\}$ is a $1$-cocycle with coefficients 
in $\C$.
(Precisely speaking $\xi_{ij}$ depends on $t$ and $c$.) 
Since $H_{zar}^1(\x [t,c],\C)=0$, we can find constants $\xi_i \in \C$
such that $\xi_i-\xi_j=\xi_{ij}$
for $i,j \in \{1,3,12\}$.
Hence $H_i-H'_i+\xi _i$ glue together and define a regular function
on $ \x [t,c] \backslash D$.
So $H_i-H'_i+\xi _i$ is a constant in the sense that it is a function of $t$
and $c$ by Corollary(4.8) for $i=1,3,12$.


\section{Deformation of the open surface $\x[t_0,c_0] \backslash \D[t_0,c_0]$
and the Hamiltonian.}
\setcounter{equation}{0}

\par
Let us fix $c=c_0$.
Then we have a family $\x_{c_0} \backslash \D_{c_0}$ of open surfaces 
$\x[t,c_0] \backslash \D[t,c_0]$ parameterized by $t$.
Namely we have 
\begin{eqnarray}
\x_{c_0} \backslash \D_{c_0} \rightarrow \Spec \C[t].
\end{eqnarray}
Since we have a flow of the second Painlev\'e equation on 
$\x_{c_0} \backslash \D_{c_0}$ so that analytically the fibration (7.1) 
is trivial.
If we take a point $t=t_0$, then 
\begin{eqnarray*}
\left(
\x_{c_0} \backslash \D_{c_0}
\right)^{an}
\simeq
\left(
\x[t_0,c_0] \backslash \D[t_0,c_0]
\right)^{an}
\times \C.
\end{eqnarray*}
Now we consider the Kodaira-Spencer map
\begin{eqnarray*}
\varphi:T_{t_0}\rightarrow
H^1
\left(
\x[t_0,c_0] \backslash \D[t_0,c_0], \Theta
\right)
\end{eqnarray*}
associated with the fibration (7.1). Here $T_{t_0}$ denotes
the tangent space of $\Spec \C[t]$ at $t=t_0$.
The open surface $\x[t_0,c_0] \backslash \D[t_0,c_0]$ is defined by
gluing together three copies $W_1$, $W_3$, $W_{12}$ of $\C^2$ by the rule
in Theorem (4.22).
We look at the Kodaira-Spencer map.
The image $\varphi \left( \left( \frac{d}{dt}\right)_{t=t_0}\right) 
\in H^1\left( \x[t_0,c_0] \backslash \D[t_0,c_0], \Theta \right)$
is given by definition by the following 1-cocycle 
$\left\{ \theta_{j,k} \right\}
\in \Gamma (W_j \cap W_k, \Theta)$ $(j,k=1,3,12)$.
Let $y_j=f_{jk}(y_k, z_k,t)$, $z_j=g_{jk}(y_k,z_k,t)$.
Then 
\begin{eqnarray*}
\theta_{j,k}
=
\left.
\frac{\partial f_{jk}}{\partial t} \frac{\partial}{\partial y_j}
+
\frac{\partial g_{jk}}{\partial t} \frac{\partial}{\partial z_j}
\right|_{t=t_0}
\in
\Gamma(W_j \cap W_k, \Theta).
\end{eqnarray*}
In our case, we have 
\begin{eqnarray*}
&&0=\theta_{1,3}=-\theta_{3,1} \in \Gamma(W_1 \cap W_3, \Theta).\\
\smallskip \nonumber\\
&&\theta_{3,12}=\frac{1}{y_{12}^2}\frac{\partial}{\partial z_3} \in 
\Gamma(W_3 \cap W_{12}, \Theta), \hspace{1cm}
-\theta_{12,3}=\theta_{3,12}
\end{eqnarray*}
and 
\begin{eqnarray*}
0=\theta_{1,12}=-\theta_{12,1}\in \Gamma(W_1 \cap W_{12}, \Theta).
\end{eqnarray*}
The 1-cocycle $\{\theta_{i,j}\}$ is cohomologous to 0.
Since we have a symplectic structure on $\x[t_0,c_0] \backslash \D[t_0,c_0]$,
we have on $\x[t_0,c_0] \backslash \D[t_0,c_0]$ an isomorphism
\begin{eqnarray}
\Theta \simeq \Omega^1.
\end{eqnarray}
The 1-cocycle $\{\omega_{ij} \}$ with coefficients in $\Omega^1$ corresponding 
$\{ \theta_{i,j} \}$ by isomorphism (7.2) is
\begin{eqnarray*}
&&\omega_{1,3}=-\omega_{3,1}=0, \hspace{1cm}
\omega_{1,12}=-\omega_{12,1}=0\\
\smallskip \nonumber\\
&&\omega_{3,12}=\frac{1}{y_{12}^2}dy_{12}
=\frac{1}{y_3^2}dy_3
\in H^0(W_3 \cap W_{12}, \Omega^1).
\end{eqnarray*}
If we consider $\omega_j=dH_j$ on $W_j$ $(j=1,3,12)$,
then $\omega_j-\omega_k=\omega_{ik}$ by Lemma (6.2).
So the cohomology class determined by $\{\omega_{ij} \} \in 
H^1\left( \x[t_0,c_0] \backslash \D[t_0,c_0], \Omega^1 \right)$ is $0$
and hence
the Kodaira-Spencer class $\{\theta_{ij}\} \in 
H^1\left( \x[t_0,c_0] \backslash \D[t_0,c_0], \Theta \right) $ is $0$.
So we can express this fact by saying that the Hamiltonian functions $H_j$
trivialize the Kodaira-Spencer class of the family (7.1).

\section{Cohomology groups.}
\setcounter{equation}{0}

\par 
We fix $t_0$, $c_0$ and denote the rational surface $\x [t_0,c_0]$ by $X$ and  $ \D [t_0, c_0]$ by 
$D$.  

\begin{Lemma}
\begin{eqnarray*}
\chi(\Theta_X)=-10.
\end{eqnarray*}
\end{Lemma}

\Proof
Since $X$ is a rational surface, $\chi({\cal O}_X)=1$.
For a vector bundle $E$ of rank 2 over $X$, the Riemann-Roch theorem tells us
\begin{eqnarray*}
\chi(E)=\frac{1}{2}
\left( c_1(E)^2-2c_2(E)\right)+\frac{1}{2}c_1(E)c_1(\Theta_X)+2\chi({\cal O}_X),\end{eqnarray*} 
and hence
\begin{eqnarray*}
\chi(E)=\frac{1}{2}
\left( c_1(E)^2-2c_2(E)\right)+\frac{1}{2}c_1(E)c_1(\Theta_X)+2.
\end{eqnarray*}
If we take $\Theta_X$ as $E$,
\begin{eqnarray*}
\chi(\Theta_X)=\frac{1}{2}
\left( c_1(\Theta_X)^2-2c_2(\Theta_X)\right)
+\frac{1}{2}c_1(\Theta_X)c_1(\Theta_X)+2
\end{eqnarray*}
so that
\begin{eqnarray*}
\chi(\Theta_X)=\frac{1}{2}
\left( c_1(\Theta_X)^2-2c_2(\Theta_X)\right)
+\frac{1}{2}c_1(\Theta_X)^2+2.
\end{eqnarray*}
Since $c_1(\Theta_X)=-K$, $c_1(\Theta_X)^2=0$ and hence
\begin{eqnarray}
\chi(\Theta_X)=-c_2(\Theta_X)+2.
\end{eqnarray}
It follows from Noether's formula
\begin{eqnarray}
1-q(X)+p_q(X)&=&\frac{1}{12}\left( c_1(\Theta_X)^2+c_2(\Theta_X)\right)
\nonumber \\
\smallskip \nonumber\\
12&=&c_2(\Theta_X).
\end{eqnarray}
Now the lemma follows from (8.2) and (8.3).
\begin{Proposition}
\begin{eqnarray*}
h^1(X, \Theta_X(- \log D))=2.
\end{eqnarray*}
\end{Proposition}

\Proof
By the Serre duality
\begin{eqnarray*}
h^1\left(X, \Theta_X(-log D) \right)=
h^1\left(X,K \otimes \Omega_X^1(\log D) \right)
\end{eqnarray*}
and hence we have to calculate 
$h^1\left(X,K \otimes \Omega_X^1(\log D) \right)$.
We have an exact sequence
\begin{eqnarray*}
0 \rightarrow \Omega_X^1 \rightarrow \Omega_X^1(\log D) \rightarrow
\bigoplus_{i=1}^8 {\cal O}_{D_i} \rightarrow 0.
\end{eqnarray*}
Tensoring $K$ with the exact sequence, we get
\begin{eqnarray*}
0 \rightarrow K \otimes \Omega_X^1 
\rightarrow K \otimes \Omega_X^1(\log D) 
\rightarrow \bigoplus_{i=1}^8 K \otimes {\cal O}_{D_i} \rightarrow 0.
\end{eqnarray*}
Since $D_i$ is a $-2$-curve, we have 
$K \otimes {\cal O}_{D_i} \simeq {\cal O}_{D_i}$
and consequently we have
\addtocounter{equation}{1}
\setcounter{meq}{\value{equation}}
\setcounter{equation}{0}
{
\renewcommand{\theequation}{\thenewsection.\themeq.\arabic{equation}}
\begin{eqnarray}
0 \rightarrow K \otimes \Omega_X^1 
\rightarrow K \otimes \Omega_X^1(\log D) 
\rightarrow \bigoplus_{i=1}^8 {\cal O}_{D_i} \rightarrow 0.
\end{eqnarray}
This exact sequence gives the long exact sequence
\begin{eqnarray}
&&0 
\rightarrow H^0\left(K \otimes \Omega_X^1 \right) 
\rightarrow H^0 \left(K \otimes \Omega_X^1(\log D)\right) 
\rightarrow H^0 \left(\bigoplus_{i=1}^8 {\cal O}_{D_i} \right)\\
\smallskip \nonumber \\
&&\rightarrow H^1 \left(K \otimes \Omega_X^1 \right) 
\rightarrow H^1 \left(K \otimes \Omega_X^1(\log D)\right) 
\rightarrow H^1 \left(\bigoplus_{i=1}^8 {\cal O}_{D_i} \right) \nonumber\\
\smallskip \nonumber \\
&&\rightarrow H^2 \left(K \otimes \Omega_X^1 \right) 
\rightarrow H^2 \left(K \otimes \Omega_X^1(\log D)\right) 
\rightarrow 0. \nonumber
\end{eqnarray}
}
\setcounter{equation}{\value{meq}}
Since $K=-\sum_{i=0}^7 n_i D_i$ with $n_i >0$ for $0 \leq i \leq 7$,
$K \otimes \Omega_X^1(\log D)$ is a subsheaf of $\Omega_X^1$ so that 
$H^0\left( K \otimes \Omega_X^1(\log D) \right) =0$.
The $D_i$'s are isomorphic to $\P^1$ and hence 
$H^0 \left(\bigoplus_{i=1}^8 {\cal O}_{D_i} \right)=\C^8$,
$H^1 \left(\bigoplus_{i=1}^8 {\cal O}_{D_i} \right)=0$.  
So it follows from the long exact sequence
\begin{eqnarray*}
0 \rightarrow \C^8
\rightarrow H^1 \left(K \otimes \Omega_X^1 \right) 
\rightarrow H^1 \left(K \otimes \Omega_X^1(\log D)\right) 
\rightarrow 0.
\end{eqnarray*}
We had to show $h^1\left(K \otimes \Omega_X^1(\log D)\right)=2$.
To this end, it suffices to show $h^1\left(K \otimes \Omega_X^1\right)=10$
by the above exact sequence.
It follows from the Serre duality 
$h^1\left(K \otimes \Omega_X^1\right)=h^1 \left( \Theta_X\right)$ so that
we have to show 
\begin{eqnarray*}
h^1(\Theta_X)=10.
\end{eqnarray*}
We notice here by the Serre duality 
$h^2 \left( \Theta_X\right) = h^0 \left(K \otimes \Omega_X^1 \right)$.
Since $-K$ is effective $K \otimes \Omega_X^1$ is a subsheaf of $\Omega_X^1$
and since $H^0\left(\Omega_X^1\right)=0$ because $X$ is rational, 
$H^0 \left( K \otimes \Omega_X^1\right)=0$.
Therefore
\begin{eqnarray}
h^2\left( \Theta_X \right)=0.
\end{eqnarray}
Now (8.4) follows from Theorem (5.26), Lemma (8.1) and (8.6).
\begin{Corollary} to the proof.
$H^0\left( \Theta_X (-\log D)\right)
=H^2\left( \Theta_X (-\log D)\right)=0.
$
\end{Corollary}

\Proof
In fact since 
$H^0(\Theta_X(-\log D)) \subset 
H^0\left( X \backslash D, \Theta_X \right)$,
$H^0\left( \Theta_X(-\log D)\right)=0$ by Theorem (5.26).
We noticed in the proof of Proposition (8.4) 
$$
H^0 \left( K \otimes \Omega_X^1(\log D) \right)=0,
$$
 implies by the
Serre duality $H^2 \left( \Theta_X(-\log D) \right)=0$.
\begin{Corollary} to the proof.
If $c_0=0$, we know that there is a $-2$-curve 
$C'_2$ on $X[t_0,0] \backslash D$.
We have 
\begin{eqnarray*}
h^1\left( X, \Theta_X \left(-\log(D+C'_2) \right) \right)=1
\end{eqnarray*}
and 
\begin{eqnarray*}
h^0\left( \Theta_X\left( - \log (D+C'_2)\right) \right)
=
h^2\left( \Theta_X\left( - \log (D+C'_2)\right) \right)
=0.
\end{eqnarray*}
\end{Corollary}

\Proof
We argue as in the proof of Proposition (8.4).
We replace $D$ by $D+C'_2$ and consider the exact sequence (8.5.1)
and the long exact sequence (8.5.2).
Then we have to show 
$H^0\left( K \otimes \Omega_X^1\left( \log(D+C'_2)\right) \right)=0$.
Since $K \otimes \Omega_X^1\left( \log(D+C'_2)\right)$ is a subsheaf
of $\Omega_X^1(\log C'_2)$,
it is sufficient to show $H^0\left( \Omega_X^1(\log C'_2)\right)=0$.
In fact we have a commutative diagram
\begin{eqnarray*}
\begin{array}{ccccccccc}
0 
&\rightarrow 
&H^0(\Omega_X^1)
&\rightarrow
&H^0\left( \Omega_X^1(\log C'_2)\right)
&\rightarrow
&H^0\left( {\cal O}_{C'_2}\right)
&\simeq
&\C\\
&&\downarrow&&\downarrow&&\parallel&&\\
&&H^1(X,\C) 
&\rightarrow&
H^1\left( X \backslash C'_2, \C \right)
&\rightarrow&
H^2\left( X; C'_2, \C \right)
&\rightarrow&
H^2(X,\C)
\end{array}
\end{eqnarray*}
of cohomology groups for the usual topology.
By the Poincar\'e duality we have $H^2(X;C'_2,\C) \simeq H_2(C'_2,\C)$
and hence the morphism
$H^2(X;C'_2,\C) \simeq H_2(C'_2, \C) \simeq \C \rightarrow H^2(X,\C)$ 
is injective.
So the morphism $H^1(X \backslash C'_2, \C) \rightarrow H^2(X;C'_2,\C)$
is trivial.
Since $H^0(\Omega_X^1)=0$, this shows $H^0(\Omega_X^1 (\log C'_2))=0$.
\section{Rational singular points.}
\setcounter{equation}{0}

\par
Let $Y$ be an affine surface defined over $\C$ and $P$ be a point.
We assume that $Y$ is normal, $Y \backslash \{D\}$ is smooth
and that $P$ is a rational double point of type $A_1$.
Let $f:\widetilde{Y} \rightarrow Y$ be the minimal resolution.
So $C:=\varphi^{-1}(P)$ is a curve isomorphic to $\P^1$ with $C^2=-2$.
$\Omega_Y^1$ is the sheaf of K\"ahler differentials on $Y$ over $\C$.
$\Theta_Y$ is the dual of $\Omega_Y^1$ so that 
$\Theta_Y={\rm Hom}_{{\cal O}_Y}(\Omega_Y^1, {\cal O}_Y)$. 
\begin{Proposition}
We have $$f_{*}\Theta_{\widetilde{Y}}(-\log C) \simeq \Theta_Y$$
and $$R^1f_{*}\Theta_{\widetilde{Y}}(-\log C)=0.$$
\end{Proposition}

\Proof
It follows form  an exact sequence
$$
0 \to \Theta_{\widetilde{Y}} (-\log C) \to \Theta_{\widetilde{Y}} \to N_C \to
0,
$$
the long exact sequence
\begin{eqnarray*}
& 0  \to  
f_* \Theta_{\widetilde{Y}} (-\log C) \to f_* \Theta_{\widetilde{Y}} \to f_* 
N_C \\
&  \to  R^1f_* \Theta_{\widetilde{Y}} (-\log C) \to 
R^1f_* \Theta_{\widetilde{Y}} 
\stackrel{i}{\longrightarrow}  
R^1f_* N_C \to 0,
\end{eqnarray*}
$N_C$ being the normal bundle of the curve $C$.
Since $N_C\simeq {\cal O}_C(-2)$, we have $f_*N_C =0.$
By Proposition(1.2) in [BW],  
$$
f_*\Theta_{\widetilde{Y}}\simeq \Theta_Y.
$$
Moreover for the minimal resolution of a rational double point, we have an isomorphism 
$H^1(\Theta_{\widetilde{Y}}) \simeq H^1(N_C)$
by (1.8) in [BW]. 
Hence the morphism $i$ in the long exact sequence above is an isomorphism
and consequently 
$$
f_*\Theta_{\widetilde{Y}}(-\log C) \simeq f_* \Theta_{\widetilde{Y}}\simeq \Theta_Y,  \qquad R^1f_*\Theta_{\widetilde{Y}}( - \log C) = 0.
$$
\section{Calculation of $\Ext$.}
\setcounter{equation}{0}

\par
On rational surface $\x[t_0,0]$,
there is a curve $C'_2$ isomorphic $\P^1$ with $\left( C'_2\right)^2=-2$
and $C'_2 \cap D =\emptyset$.
We can contract the curve $C'_2$ on $\x[t_0,0]$ to 
a rational singular point $P$ to get a rational surface $\Y[t_0,0]$ 
with a rational double point $P$ of
type $A_1$.
We have a morphism
\begin{eqnarray*}
f:\x[t_0,0] \rightarrow \Y[t_0,0],
\end{eqnarray*}
which is the minimal resolution of $\Y[t_0,0]$.
Since the curve $D$ does not intersect $C'_2$, $f$ is an isomorphism on a 
neighborhood of $D$.
So we denote the image $f(D)$ on $\Y[t_0,0]$ again by $D$.
We compare the Kuranishi family of the pair $\left( \x[t_0,0], D \right)$
and that of the pair $\left( \Y[t_0,0], D \right)$.
We denote $\x[t_0,0]$ by $X$ and $\Y[t_0,0]$ by $Y$.
\begin{Lemma}
\begin{eqnarray*}
H^2(Y, \Theta_Y(-\log D))=0
\end{eqnarray*}
\end{Lemma}

\Proof
It follows from Proposition (9.1) that the spectral sequence
\begin{eqnarray*}
E_2^{p,q}
&=&H^p\left(Y, R^qf_{*}\Theta_X\left(-\log(D+C'_2)\right) \right)\\
\smallskip \nonumber \\
&\Rightarrow& 
H^n\left(X, \Theta_X \left(-\log(D+C'_2)\right) \right) \nonumber
\end{eqnarray*}
degenerates and that $f_{*}\Theta_X \left( -\log(D+C'_2) \right) 
=\Theta_Y(-\log D)$.
So we have
\begin{eqnarray*}
H^p\left(Y, \Theta_Y(-\log D)\right)
=H^p\left(X, \Theta_X\left(-\log(D+C'_2)\right)\right)
\end{eqnarray*}
for every integer $p \geq 0$.
The lemma now follows from Corollary (8.8).
\begin{Proposition}
The spectral sequence
\begin{eqnarray*}
E_2^{p,q}=H^p\left(Y, 
{\cal E}{xt}_{{\cal O}_Y}^q
\left(
\Omega_Y^1(\log D), {\cal O}_Y 
\right)
\right)\\
\smallskip \nonumber\\
\Rightarrow
\Ext_{{\cal O}_Y}^{p+q} \left( \Omega_Y^1(\log D), {\cal O}_Y \right) \nonumber
\end{eqnarray*}
degenerates at the $E_2$-terms and we have the exact sequence
\begin{eqnarray*}
0 \rightarrow H^1\left( Y, \Omega_Y^1(\log D)\right)
\rightarrow \Ext_{{\cal O}_Y}^1\left( \Omega^1(\log D), {\cal O}_Y\right)
\rightarrow \Ext_{{\cal O}_{Y,P}}^1
\left( \Omega_{Y,P}^1, {\cal O}_{Y,P}\right)
\rightarrow 0.
\end{eqnarray*} 
\end{Proposition}

\Proof
Since $Y$ is locally a complete intersection,
the projective dimension of $\Omega_Y^1(-\log D)$ is equal to
the projective dimension of $\Omega_Y^1 \leq 1$ so that
$$
{\cal E}{xt}_{{\cal O}_Y}^q \left( \Omega_Y^1(-\log D), {\cal O}_Y\right)=0
$$
for $q \geq 2$. 
${\cal E}{xt}_{{\cal O}_Y}^1 \left(\Omega_Y^1(-\log D), {\cal O}_Y \right)$
is supported on the singular point $P$.
Lemma (10.1) shows $E_2^{0,2}=0$.
So the spectral sequence degenerates at the $E_2$-terms.
Thus $\Ext_{{\cal O}_Y}^2\left( \Omega^1(\log D), {\cal O}_Y \right)=0$
and we have an exact sequence
\begin{eqnarray*}
0 \rightarrow E_2^{1,0} \rightarrow 
\Ext_{{\cal O}_Y}^1 \left( \Omega_Y^1(-\log D), {\cal O}_Y \right)
\rightarrow E_2^{0,1} \rightarrow 0,
\end{eqnarray*}
which is nothing but the exact sequence of the proposition.

\section{Family $\x[t,c]$ of rational surfaces.}

\setcounter{equation}{0}
\begin{Theorem}
The family $(\x[t,c], \D[t,c]) \rightarrow \Spec \C[t,c]$
of pairs a surface and divisor
is a Kuranishi family at every point 
$(t_0,c_0) \in \left( \C[t,c]\right)^{an} = \C^2$. 
\end{Theorem}

\Proof Again we denote $\D[t_0, c_0]$ by $D$.  
The theorem says that the Kodaira-Spencer map
\begin{eqnarray}
\rho:T_{\C^2,P}\rightarrow H^1\left(\x[t_0,c_0], \Theta(-\log D) \right)
\end{eqnarray}
is an isomorphism of vector spaces,
where $T_{\C^2,P}$ is the tangent space of 
$\left( \Spec \C[t,c] \right)^{an}=\C^2$ at a point $P=(t_0,c_0)$.
It follows from Proposition (8.4) that we have to show that the
image of Kodaira-Spencer map is of dimension 2.
As we have show in \S2, \cite{Utr}, 
the surface $\x[t_0,c_0]$ is covered by 12 affine open sets $W_i$ isomorphic to
$\A^2$ $(1 \leq i \leq 12)$.
The image
\begin{eqnarray*}
\rho \left( \left. \left( \frac{\partial}{\partial c}\right)\right|_{(t,c)
=(t_0,c_0)}\right)
\in
H^1\left(\x[t_0,c_0], \Theta(-\log D) \right)
\end{eqnarray*}
is represented by a \v{C}ech 1-cocycle $a=\{a_{i,j}\}_{1 \leq i,j \leq 12}$
with coefficients in $\Theta(- \log D)$ with respect to the covering
$X[t_0,c_0]=\cup_{i=1}^{12} W_i$ so that
\begin{eqnarray*}
a_{i,j} \in H^0 \left( W_i \cap W_j, \Theta(-\log D)\right).
\end{eqnarray*}
Similarly, the image
\begin{eqnarray*}
\rho \left( \left. \left( \frac{\partial}{\partial t}\right)\right|_{(t,c)
=(t_0,c_0)}\right)
\in
H^1\left(\x[t_0,c_0], \Theta(-\log D) \right)
\end{eqnarray*}
is given by \v{C}ech 1-cocycle 
\begin{eqnarray*}
b=\{b_{i,j}\}_{1 \leq i,j \leq 12}
\end{eqnarray*}
with
\begin{eqnarray*}
b_{i,j} \in H^0 \left( W_i \cap W_j, \Theta(-\log D)\right).
\end{eqnarray*}
We assume that a linear combination $\lambda a+ \mu b$, which is a 
\v{C}ech 1-cocycle with coefficients in $\Theta(-\log D)$,
is cohomologous to $0$ for complex numbers $\lambda$, $\mu$ and we show that
$\lambda=\mu=0$.
We recall that $\x[t_0,c_0] \backslash D$ is covered by $W_1$, $W_3$ 
and $W_{12}$.
We have a canonical map
\begin{eqnarray*}
H^1\left( \x[t_0,c_0], \Theta(-\log D) \right)
\rightarrow
H^1\left( \x[t_0,c_0] \backslash D, \Theta(-\log D) \right)
\end{eqnarray*}
The image of a \v{C}ech 1-cocycle $f=\{f_{i,j}\}$ under this map is denoted
by $\bar{f}=\{\bar{f}_{i,j}\}$.
Since as we have seen in \S7, $\bar{b}$ is cohomologous to $0$,
it follows from the assumption that 
$\lambda \bar{a}=\lambda \bar{a} +\mu \bar{b}$ is cohomologous to $0$.
By the coordinate transformation between $W_1$ and $W_3$ given in \S5,
we get 
\begin{eqnarray*}
\frac{\partial z_1}{\partial c}=
\frac{\partial}{\partial c}(y_3c-y_3z_3)=y_3=\frac{1}{y_1}
\end{eqnarray*}
and $\partial y_1 / \partial c =0$.
So by the definition of the Kodaira-Spencer map, we have 
\begin{eqnarray}
\bar{a}_{1,3}=\frac{1}{y_1}\frac{\partial}{ \partial z_1}\in H^0
\left(W_1 \cap W_3,\Theta(-\log D)\right).
\end{eqnarray} 
Similarly
\begin{eqnarray}
\bar{a}_{1,12}=-\frac{1}{y_1}\frac{\partial}{ \partial z_1}\in H^0
\left(W_1 \cap W_{12}, \Theta(-\log D)\right).
\end{eqnarray}
On the other hand, we have the symplectic structure on 
$\x[t_0,c_0] \backslash D=W_1 \cup W_2 \cup W_3$ given by
$dy_1 \wedge dz_1=dy_3 \wedge dz_3=dy_{12} \wedge dz_{12}$ so that
we have an isomorphism
\begin{eqnarray}
\varphi:\Theta \simeq \Omega^1
\end{eqnarray}
on $\x[t_0,c_0] \backslash D=W_1 \cup W_2 \cup W_3$.
Under this isomorphism $\partial/\partial y_i$ corresponds to $dz_i$
and $\partial/\partial z_i$ to $-dy_i$.
We set 
$\bar{\alpha}_{i,j}:=\varphi(\bar{a}_{i,j}) \in H^0(W_i \cap W_j, \Omega^1)$
for $i,j=1,3,12$.
In fact we have explicitly
\begin{eqnarray}
\bar{\alpha}_{1,3}=-\frac{1}{y_1}dy_1, \hspace{1cm}
\bar{\alpha}_{1,12}=\frac{1}{y_1}dy_1.
\end{eqnarray}
Since $\lambda \bar{\alpha}$ is cohomologous to $0$,
there exist $\bar{\alpha}_i \in H^0(W_i, \Omega^1)$ for $i=1,3,12$
such that
\begin{eqnarray}
\bar{\alpha}_1-\bar{\alpha}_3=-\frac{\lambda}{y_1}dy_1,
\hspace{1cm}
\bar{\alpha}_1-\bar{\alpha}_{12}=\frac{\lambda}{y_1}dy_1.
\end{eqnarray}
In view of (11.6) and (11.7), 
$d\bar{\alpha}_1=d\bar{\alpha}_3=d\bar{\alpha}_{12}$
glue together and define a $2$-form on $\x[t_0,c_0] \backslash D$.
So by Corollary (4.7), $d\bar{\alpha}_1=d\bar{\alpha}_3=d\bar{\alpha}_{12}
=\nu \omega$ for a complex number $\nu$.
Since $\omega$ is not exact by Corollary (5.26), we have $\nu=0$ and 
$d\bar{\alpha}_1=d\bar{\alpha}_3=d\bar{\alpha}_{12}=0$.
Now since $W_i \simeq \A^2$, by a theorem of Grothendieck,
there exist polynomials 
$f_1 \in \C[y_1,z_1]$, 
$f_3 \in \C[y_3,z_3]$, 
$f_{12} \in \C[y_{12},z_{12}]$ 
such that $df_i=\bar{\alpha}_i$ for $i=1,3,12$.
So it follows from (11.7)
\begin{eqnarray*}
df_1-df_2=-\frac{\lambda}{y_1}dy_1
\end{eqnarray*}
on $W_1 \cap W_3=\Spec \C[y_1,1/y_1,z_1]$.
Since $-\frac{\lambda}{y_1}dy_1$ is exact on $W_1 \cap W_3$ if and only if 
$\lambda=0$.
We conclude $\lambda=0$.
Now it remains to show $\mu=0$.
Let us  consider the canonical map
\begin{eqnarray*}
H^1\left(\x[t_0,c_0], \Theta(-\log D) \right)
\rightarrow
H^1\left( W_1 \cup W_2 \cup W_3 \cup W_{12}, \Theta(-\log D)\right).
\end{eqnarray*}
The image of a 1-\v{C}ech cocycle $f=\{f_{i,j}\}$ under this map will be 
denoted by $\widetilde{f}=\{\widetilde{f}_{i,j}\}$.
We calculated $\widetilde{b}$ for $i=1,3,12$ in \S6.
We recall the coordinate transformation between $W_1$ and $W_2$,
which are a part of coverings of the starting ruled surface in the 
construction of $X[t_0,c_0]$, is given by $y_1=y_2$ and $z_1=1/z_2$
(cf. \S1).
It follows from the definition of the Kodaira-Spencer map 
$\widetilde{b}_{1,2}=0$.
Since $\mu b$ is cohomologous to $0$, 
$\mu \widetilde{b}$ is cohomologous to $0$ too.
We can find $\widetilde{b}_i \in H^0\left(W_i,\Theta(-\log D) \right)$
for $i=1,2,3,12$ such that
\begin{eqnarray}
\widetilde{b}_i-\widetilde{b}_j=\mu \widetilde{b}_{i,j}
\end{eqnarray}
for $i,j=1,2,3,12$.
Since $H^0\left( \x[t_0,c_0] \backslash D, \Theta \right)
=H^0\left( \x[t_0,c_0] \backslash D, \Theta(-\log D)\right)=0$,
$\widetilde{b}_i \in H^0\left( W_i, \Theta \right)$ is uniquely determined for
$i=1,3,12$.
Namely we have 
\begin{eqnarray*}
\widetilde{b}_{t,i}
=\mu \left( 
\frac{\partial H_i}{\partial z_i}
\frac{\partial}{\partial y_i}
-
\frac{\partial H_i}{\partial y_i}
\frac{\partial}{\partial z_i}
\right)
\end{eqnarray*}
for $i=1,3,12$, where $H_i$ is the Hamiltonian as we have seen in \S7.
In particular
\begin{eqnarray*}
\widetilde{b}_1=\mu
\left\{
\left(y_1^2+z_1+\frac{t}{2} \right)\frac{\partial}{\partial y_1}
+(c-2y_1z_1)\frac{\partial}{\partial z_1}
\right\}.
\end{eqnarray*}
Since as we noticed above $\widetilde{b}_{1,2}=0$, we have
\begin{eqnarray*}
\widetilde{b}_1-\widetilde{b}_2=0,
\end{eqnarray*}
which means $\widetilde{b}$ is regular on $W_2$. We have
\begin{eqnarray*}
\widetilde{b}_1=
\mu \left\{
\left(y_2^2+\frac{1}{z_2}+\frac{t}{2}\right)
\frac{\partial}{\partial y_2}
+\left( c-2y_2\frac{1}{z_2}\right)
\left( -z_2^2\right)
\frac{\partial}{\partial z_2}
\right\}
\end{eqnarray*}
on $W_2$.
So if $\widetilde{b}_1$ is regular on $W_2$, then $\mu=0$.
This is what we had to prove.
\par
It is an interesting problem to give a pair of algebraic varieties 
$(U,\, V)$
defined over $\C$
such that the algebraic varieties $U$ and $V$ are not isomorphic one 
another but the associated analytic spaces $U^{an}$ and $V^{an}$ are
isomorphic.
A well-known example due to Serre is related with algebraic groups(See [H], Chapter VI, \S 3). 
The family $(\x[t,c], \, \D[t,\, c])$ provides us many such examples.
Let us choose a point $(t_0,\, c_0)\in \C^2.$ 
Let $t_1\not= t_0$ be an arbitrary point in a small neighbourhood 
of $t_0\in \C.$
The pair 
$\x[t_0 , \, c_0] \setminus  \D [ t_0 ,\, c_0], \, \x[t_1, \, c_0]
\setminus
  \D[ t_1,\, c_0])$
is a such pair.
In fact, $\x[t_0, \, c_0]\setminus  \D[t_0,\, c_0]$
is not isomorphic to 
$\x[t_1, \, c_0]\setminus  \D[t_1,\, c_0]$
as algebraic varieties.
In fact if we had an isomorphism 
$$
\varphi :\x[t_0, \, c_0]\setminus  \D[t_0,\, c_0] \to \x[t_1, \, c_0]\setminus  \D[t_1,\, c_0]
$$
as algebraic varieties, then 
since 
the complementary divisors $\D[t_0, \, c_0]$ and $\D[t_1, \, c_0]$ 
consist of $-2$-curves,
we could extend $\varphi$ to an isomorphism
$$
\bar{\varphi}:\x[t_0, \, c_0] \to \x[t_1, \, c_0]  
$$
of algebraic varieties
so that the pair 
$
(\x[t_0, \, c_0],\,  \D[t_0,\, c_0])
$
is isomorphic to the pair 
$
(\x[t_1, \, c_0], \,  \D[t_1,\, c_0]).
$
This contradicts Theorem(11.1).
Integration of the second Painlev\'e equation that has no movable singular points gives an analytic isomorphism
$$
\x[t_0, \, c_0]\setminus  \D[t_0,\, c_0] \simeq \x[t_1, \, c_0]\setminus  \D[t_1,\, c_0].
$$

\end{document}